\documentclass[11pt]{article}
\usepackage{amsmath}
\usepackage{amssymb}
\usepackage{color}
\usepackage{latexsym}
\usepackage{comment}

\newtheorem{assumption}{Assumption}
\def\qed{ \ \vrule width.2cm height.2cm depth0cm\smallskip}
\newenvironment{proof}{\noindent {\bf Proof.\/}}{$\qed$\vskip 0.1in}

\newcommand{\la}{\langle}
\newcommand{\ra}{\rangle}

\newcommand{\ba}{\begin{array}}
\newcommand{\ea}{\end{array}}
\newcommand{\be}{\begin{equation}}
\newcommand{\ee}{\end{equation}}
\newcommand{\bea}{\begin{eqnarray}}
\newcommand{\eea}{\end{eqnarray}}
\newcommand{\beaa}{\begin{eqnarray*}}
\newcommand{\eeaa}{\end{eqnarray*}}

\def\dbC{\mathbb{C}}

\def\dbE{\mathbb{E}}
\def\dbF{\mathbb{F}}

\def\dbK{\mathbb{K}}
\def\dbL{\mathbb{L}}

\def\dbN{\mathbb{N}}
\def\dbP{\mathbb{P}}
\def\dbR{\mathbb{R}}
\def\dbS{\mathbb{S}}

\def\db0{\bf{0}}

\def\a{\alpha}
\def\b{\beta}
\def\g{\gamma}
\def\d{\delta}
\def\e{\varepsilon}

\def\l{\lambda}

\def\si{\sigma}
\def\t{\tau}
\def\f{\varphi}

\def\o{\omega}

\def\G{\Gamma}

\def\L{\Lambda}

\def\O{\Omega}
\def\cA{{\cal A}}

\def\cD{{\cal D}}
\def\cE{{\cal E}}

\def\cH{{\cal H}}

\def\cJ{{\cal J}}
\def\cK{{\cal K}}
\def\cL{{\cal L}}
\def\cM{{\cal M}}

\def\cO{{\cal O}}
\def\cP{{\cal P}}

\def\cS{{\cal S}}
\def\cT{{\cal T}}
\def\cU{{\cal U}}

\def\cY{{\cal Y}}
\def\cZ{{\cal Z}}
\def\ch{\textsc{h}}

\def\q{\quad}
\def\qq{\qquad}

\def\pa{\partial}
\def\cd{\cdot}

\def\tr{\hbox{\rm tr}}
\def\qed{ \hfill \vrule width.25cm height.25cm depth0cm\smallskip}

\newcommand{\basa}{\begin{assumption}}
\newcommand{\easa}{\end{assumption}}

\newcommand{\bas}{\begin{assum}}
\newcommand{\eas}{\end{assum}}

\def\pa{\partial}

 \def\cd{\cdot}

\def\tr{\hbox{\rm tr$\,$}}

\def\wh{\widehat}

\def\1{{\bf 1}}

\def\:{\!:\!}
\def\reff#1{{\rm(\ref{#1})}}
\def \proof{{\noindent \bf Proof\quad}}

\def \Usup{\overline{\cU}}
\def \Usub{\underline{\cU}}
  
at 9pt
\begin{document}
\newtheorem{thm}{Theorem}[section]
\newtheorem{lem}[thm]{Lemma}
\newtheorem{cor}[thm]{Corollary}
\newtheorem{prop}[thm]{Proposition}
\newtheorem{rem}[thm]{Remark}
\newtheorem{eg}[thm]{Example}
\newtheorem{defn}[thm]{Definition}
\newtheorem{assum}[thm]{Assumption}

\renewcommand {\theequation}{\arabic{section}.\arabic{equation}}
\def\thesection{\arabic{section}}

\title{\bf Viscosity solutions of obstacle problems for Fully nonlinear path-dependent PDEs}

\author{Ibrahim  {\sc Ekren}\footnote{The author would like to thank Nizar Touzi and Jianfeng Zhang for helpful discussions.}   \footnote{University of Southern California, Department of Mathematics, ekren@usc.edu}   
}\maketitle

\begin{abstract}
In this article, we adapt the definition of viscosity solutions to the obstacle problem for fully nonlinear path-dependent PDEs with data uniformly continuous in $(t,\o)$, and generator Lipschitz continuous in $(y,z,\g)$. We prove that our definition of viscosity solutions is consistent with the classical solutions,  and satisfy a stability result. We show that the value functional defined via the second order reflected backward stochastic differential equation is the unique viscosity solution of the variational inequalities. 
\end{abstract}
 


\vfill\eject
\section{Introduction}
Since the seminal work of Pardoux and Peng \cite{BSDE}, backward stochastic differential equations(BSDEs) have found many areas of application. In \cite{EKKPQ}, El Karoui et al. introduced a new kind of BSDE, called reflected BSDEs, where the solution is forced to stay above a barrier. In the Markovian framework, they have proven that the value function defined through the RBSDE is the unique viscosity solution of an obstacle problem for a semi-linear parabolic partial differential equation, hence extending the well-known Feynman-Kac formula to the associated variational inequalities. 

In his recent work  \cite{DPath}, Dupire gives the definition of derivatives on the path space and proves a functional It\^o's formula. Using those derivatives, in  \cite{DPath}, and \cite{CF}, the authors derive and study functional differential equations, which extends the Feynman-Kac formula to a non-Markovian case. In \cite{PENG}, Peng  proposes the notion of path-dependent partial differential equations (PPDEs) in nonlinear framework.  In \cite{EKTZ}, \cite{ETZ1}, \cite{ETZ2}, Ekren et al. proposed a definition for viscosity solution of PPDEs. 

Our objective in this paper is to adapt the definition of viscosity solutions of PPDEs given in \cite{ETZ1} to an obstacle problem for a fully nonlinear PPDE, for which we give assumption under which wellposedness holds. In order to achieve our objective, especially to tackle with the lack of local compactness of the space of paths, we will use similar ideas as in \cite{ETZ2}. In our case, the main difficulty is to produce a sequence of "smooth" subsolutions and supersolutions of the obstacle problem that converge to the value functional, however in general the solutions of obstacle problem of PDEs do not have $C^{1,2}$ regularity. To overcome this difficulty, we use a penalization approach and a change of variable which allows us to have "smooth" solution to the obstacle problem.  Another main difference with \cite{ETZ2} is the fact that the PPDE studied in this paper has a stochastic representation with a second order reflected BSDE. This will allow us to prove the regularity of the functional without requiring an assumption similar to Assumption 3.5 of \cite{ETZ2}. 

The paper is organized as follows. In section 2, we introduce the main notations and assumptions that we will use. In section 3 we  introduce our functional of interest as the supremum of solutions of RBSDEs and give some regularity results for this functional. In section 4, we introduce the PPDE we want to study and give the definition of viscosity solution for this kind of PPDEs and prove some preliminary results on viscosity solutions. We then prove that our value functional of interest is a viscosity solution of this PPDE. 
Starting from section 5, we treat the wellposedness of the PPDE. We first prove our partial comparison result at section 5. In section 6, a stability theorem is proven. In section 7 we prove, using a modification of Perron's approach,  the general comparison result (without requiring any smoothness  of the sub and supersolutions).
\section{Notations}

We fix $T>0$, the time maturity, and an integer $d>0$. We denote by $\dbS^d$ the set of symmetric $d$-dimensional square matrices. For $x\in \dbR^d$, $|x|$ is the norm of $x$, and for $A,B\in \dbS^d$, $ A:B:=trace(AB)$. For any matrix $M$, $M^*$ denotes its transpose. We work on the canonical space $\O:=\{\o\in C([0,T],\dbR^d):\o_0=0\}$ of $d$-dimensional continuous paths. $B$ denotes the canonical process on this space, $\dbF=\{\dbF_s\}_{{s\in[0,T]}}$ is the filtration generated by $B$, and $\dbP_0$ is the Wiener measure. For $\o\in \O$ and $t\in[0,T]$, the stopped path $\o_{.\wedge t}\in\O$ is defined as follows : 
\beaa
&&\o_{.\wedge t} (s)=\o_s, \mbox { for } 0\leq s \leq t,\\
&&\o_{.\wedge t} (s)=\o_t, \mbox { for } t\leq s \leq T.
\eeaa
We denote $\L:=\{(t,\o_{.\wedge t}):0\leq t\leq T,$ $ \o\in\O\}$. In the sequel, we will denote a generic element of $\L$ as $(t,\o)\in[0,T]\times \O$, this notation means that we ignore the values of $\o$ after $t$ and identify $(t,\o')\in [0,T]\times \O$ with $(t,\o)\in[0,T]\times \O$ if $\o_{.\wedge t}=\o'_{.\wedge t}$. We define the following $||.||_T$ and $d_\infty$ metrics on respectively $\O$ and $\L$: 
\beaa
&& || \o ||_T:=\sup_{s\in [0,T]} |\o_s|,\quad \mbox{ for } \o\in \O,\\
&& d_\infty ((t,\o),(t',\o'))=|t-t'|+||\o_{.\wedge t}-\o'_{.\wedge t'}||_T,\quad \mbox{ for } (t,\o),(t',\o')\in \L,
\eeaa
then $(\O,||.||_T)$ and $(\L,d_\infty)$ are complete metric spaces. $\dbL^0 (\dbF_T,\dbK)$ and $\dbL^0 (\L,\dbK)$(where $\dbK=\dbR,\dbR^d$ or $\dbS^d$) denote respectively the space of $\dbF_T$ measurable $\dbK$-valued random variables and $\dbF-$progressively measurable $\dbK$-valued processes. When $\dbK=\dbR$, we omit the symbol $\dbR$.

\subsection{Shifted Spaces}
For fixed $s\leq t\in[0,T]$, we define the following shifted objects :\\
$\bullet$ $\O^t:=\{\o\in C([t,T],\dbR^d):\o_t=0\}$.\\
$\bullet$ $B^t$ is the canonical process on $\O^t$.\\
$\bullet$ $\dbF^t=\{\dbF^t_s\}_{s\in[t,T]}$ is the filtration generated by $B^t$.\\
$\bullet$ $\dbP^{t}_0$ is the Wiener measure on $\O^t$, $\dbE^t_0$ is the expectation under $\dbP^t_0$. \\
$\bullet$ We define similarly $\L^t$, $||.||^t_T$, $d^t_\infty$, $\dbL^0(\dbF^t_T)$ etc. In these definitions, the superscripts will generally stand for the shifted space (i.e. the beginning of times) and subscripts for the final time related to the notation.\\
$\bullet$ For $\o\in \O^s$ and $ \o'\in\O^t$, we define $\o\otimes_t\o'\in\O^s$ the concatenation of $\o$ and $\o'$ at $t$ by: 
\beaa
(\o\otimes_t \o') (r) := \o_r\1_{[s,t)}(r) + (\o_{t} + \o'_r)\1_{[t, T]}(r),
&\mbox{for all}&
r\in [s,T].
\eeaa
$\bullet$ For $\xi\in\dbL^0(\dbF^{s}_{T})$, $X\in\dbL^0 (\L^s)$, and a fixed path $\o\in\O^s$,we define the shifted random variable $\xi^{t,\o}\in\dbL^0 (\dbF^{t}_{T})$  and   process $X^{t,\o}\in\dbL^0 (\L^t)$ by:
\beaa
\xi^{t, \o}(\o') :=\xi(\o\otimes_t \o'), \q X^{t, \o}(\o') := X(\o\otimes_t \o'),
&\mbox{for all}&
\o'\in\O^t.
\eeaa
$\bullet$ Finally, we denote by $\cT$ the set of $\dbF$-stopping times, $\cT_+$ the set of positive $\dbF$-stopping times, and $\cH\subset\cT$ the subset of those hitting times $\ch$ of the form
\bea
\label{cT}
\ch := \inf\{t: B_t \in O^c\} \wedge t_0 = \inf\{ t: d(\o_t, O^c) = 0\} \wedge t_0,
\eea
for some $0< t_0\le T$, and some open and convex set $O  \subset \dbR^d$ containing ${\bf 0}$ with $O^c := \dbR^d\setminus O$.
\beaa
\ch>0, \mbox{$\ch$ is lower semi-continuous, and $\ch_1 \wedge \ch_2 \in \cH$ for any $\ch_1, \ch_2\in\cH$}.
\eeaa
$\cT^t$ and $\cH^t$ are defined in the same spirit. For any $\t\in \cT$ (resp. $\ch\in\cH$) and any $(t,\o)\in\L$ such that $t<\t(\o)$ (resp. $t<\ch(\o)$), it is clear that $\t^{t,\o}\in \cT^t$ (resp. $\ch^{t,\o}\in\cH^t$). For $t\in[0,T]$ and $\d>0$ we define the hitting times $\ch^t_\d\in \cH^t$, that we will use several times, as follows : 
\bea\label{definition_ch}
\ch^t_\d:=\inf\{s\geq t :||B^t||_s=\d\}\wedge (t+\d)\wedge T.
\eea
Notice that $O$ is open and contains ${\bf 0}$, and $t_0>t$. Therefore, for all $\ch\in\cH^t$  there exist $\d>0$ such that 
\bea\label{che}
t< \ch^t_\d\leq \ch. 
\eea
The class $\cH$ and especially, the stopping time $\ch^t_\d$ will be our mains tools for studying processes locally to the right in time in the space $\L$.

We shall use the following type of regularity, which is stronger than the right continuity of a process in a standard stochastic analysis sense.
\begin{defn}
\label{defn-rightcont}
We say  a process $u \in \dbL^0 (\L)$ is right continuous under $d_\infty$ if for any $(t,\o) \in \L$ and any $\e>0$, there exists $\d>0$ such that, for any $(s, \tilde\o)\in \L^t$ satisfying $d^t_\infty((s,  \tilde\o), (t, {\bf 0})) \le \d$, we have $|u^{t,\o}(s,\tilde \o) - u(t,\o)|\le \e$.
\end{defn} 

We now define the following sets of functionals which are the equivalents of semi continuous and continuous functions in the viscosity solutions theory of PDEs. Notice that for a mapping $u:\L\rightarrow \dbK$, $\dbF$-progressive measurability implies that $u(t,\o)=u(t,\o_{.\wedge t})$ for all $(t,\o)\in [0,T]\times \O$.

\begin {defn}
(i) $\underline \cU\subset \dbL^0(\L)$ is the set of processes $u$ that are bounded from above, right continuous under $d_\infty$, and such that there exist a modulus of continuity $\rho$ verifying for any $(t,\o), (t',\o')\in\L$:
 \bea\label{USC}
 u(t,\o) - u(t',\o') 
 \le 
 \rho\Big(d_\infty\big((t,\o), (t',\o')\big)\Big)
 ~\mbox{whenever}~t\le t'.
 \eea 
(ii) $\overline \cU\subset \dbL^0(\L)$ is the set of processes $u$ such that $-u\in\underline \cU$. \\
(iii) $C^0(\L,\dbK)$(respectively, $C^0_b(\L,\dbK)$, $UC_b(\L,\dbK)$) is the set of $\dbF$-progressively measurable processes with values in $\dbK$ that are continuous(respectively, continuous and bounded, uniformly continuous and bounded) in $(t,\o)$ under the $d_\infty$ metric. When $\dbK=\dbR$, we simply write these sets as $C^0(\L),C^0_b(\L)$ and $UC^0_b(\L)$. 
\end{defn}

\begin{rem}
\label{rem-C0b}
{\rm 
It is clear that $\Usub \cap \Usup = UC_b(\L)$. We also recall from \cite{ETZ0} Remark 3.2 the condition \reff{USC} implies that $u$ has left-limits and $u_{t-} \le u_t$ for all $t\in(0,T]$. 
\qed}
\end{rem}

\begin{rem}
{\rm The inequality \reff{USC} is needed to apply the results in \cite{ETZ0}. More powerful results then the one in \cite{ETZ0} are available in the literature if one wants only to study the obstacle problem for semi-linear PPDEs(see. \cite{HL}). In this particular case, it is possible to prove the comparison theorem if we define $\underline \cU$ as the class of cadlag process that are left upper semi-continuous. Notice that this last definition would not require regularity in $\o$, hence the comparison result would be more powerful than the one proven in Section 6.    
\qed}
\end{rem}

We define $\overline \cU^t$, $\underline \cU^t$, $C^0(\L^t,\dbK)$, $C^0_b(\L^t,\dbK)$ and $UC_b(\L^t,\dbK)$ in the obvious way.
It is clear that, for any $(t,\o) \in \L$ and any $u\in C^0(\L)$, we have $u^{t,\o} \in C^0(\L^t)$. The other spaces introduced before enjoy the same property. Notice also that  for $u\in C^0(\L,\dbK)$, the sample paths of $\{u(t,B)\}_{\{t\in[0,T]\}}$ are continuous.

\subsection{Nonlinear expectation}
We now a give a quick description of the the probability sets and associated capacities that we will need to define viscosity solutions of PPDEs. These sets are the ones used in \cite{ETZ1}.

For every constant $L>0$, we denote by $\cP_L$ the collection of all continuous semimartingale measures $\dbP$ on $\O$ whose drift and diffusion characteristics are bounded by $L$ and $\sqrt{2L}$, respectively.  To be precise, let $\tilde\O:=\O^3$ be an enlarged canonical space, $\tilde B:= (B, A, M)$ be the canonical processes, and $\tilde\o = (\o, a, m)\in \tilde\O$ be the paths. For any $\dbP\in \cP_L$, there exists an extension $Q$ on $\tilde\O$ such that: 
\bea
\label{cPL}
\left.\ba{c}
B=A+M,\q \mbox{ $A$ is absolutely continuous, $M$ is a martingale},\\
|\alpha^\dbP|\le L,~~{\frac12}\tr((\beta^\dbP)^2)\le L,\q \mbox{where}~ \a^\dbP_t:= \frac{d A_t}{ dt},~ \b^\dbP_t := \frac{\sqrt{d\la M\ra_t}}{dt},
\ea\right. Q\mbox{-a.s.}
\eea
Similarly, for any $t\in [0, T)$, we may define $\cP^t_L$ on $\O^t$.

We denote by $\dbL^1(\dbF^t_T,\cP^t_L)$ the set of $\xi\in \dbL^0(\dbF^t_T)$ satisfying $\sup_{\dbP\in\cP^t_L}\dbE^{\dbP}[|\xi|]<\infty$. The following nonlinear expectation will play a crucial role:
 \bea\label{cE}
 \overline{\cE}^L_t[\xi]
:=
 \sup_{\dbP\in\cP^t_L}\dbE^{\dbP}[\xi]
 ~~\mbox{and}~~
 \underline{\cE}^L_t[\xi]
 :=
 \inf_{\dbP\in\cP^t_L}\dbE^{\dbP}[\xi]
 =
 -\overline{\cE}^L_t[-\xi]
 ~~\mbox{for all}~~
 \xi\in\dbL^1(\dbF^t_T,\cP^t_L).
 \eea
We recall the following lemma whose proof can be found in \cite{ETZ1}.
\begin{lem}
\label{lem-che}
For any $\ch\in \cH$ and any $L>0$, we have $\underline \cE^L_0[\ch] >0$.
\end{lem}

\begin{defn}
Let $X\in \dbL^0(\L)$ such that $X_t\in\dbL^1(\dbF_t,\cP_L)$ for all $0\le t\le T$. We say that $X$ is an $\overline{\cE}^L-$supermartingale (resp. submartingale, martingale) if, for any $(t,\o)\in\L$ and any $\t\in \cT^t$,  $\overline{\cE}^L_t[X^{t,\o}_\t]\le$ (resp. $\ge,=$) $X_t(\o)$.
\end{defn}

We now state an important result for our subsequent analysis. Given a bounded process $X\in \dbL^0(\L)$, %$\dbF$-progressively measurable process $X$, 
consider the nonlinear optimal stopping problem 
 \bea\label{cS}
 \overline{\cS}^L_t[X](\o) 
 := 
 \sup_{\t\in \cT^t}
 \overline\cE^L_t\big[X^{t,\o}_{\t}\big]
 &\mbox{and}&
 \underline{\cS}^L_t[X](\o) 
 := 
 \inf_{\t\in \cT^t}
 \underline\cE^L_t\big[X^{t,\o}_{\t}\big], 
 ~~(t,\o)\in\L.
 \eea
By definition, we have $\overline{\cS}^L[X]\ge X$ and $\overline{\cS}^L_T[X] = X_T$. The following nonlinear Snell envelope characterization of the optimal stopping time is proven in \cite{ETZ0}.

\begin{thm}\label{thm-optimal}
Let $X\in\Usub$ be bounded, $\ch\in \cH$, and set $\wh X_t:=X_t\1_{\{t<\ch\}}+X_{\ch-}\1_{\{t\ge \ch\}}$. Define
 \beaa
 Y
 := 
 \overline\cS^L\big[\wh X\big]   
 &\mbox{and}& 
 \t^*:=\inf\{t\ge 0:Y_t=\wh X_t\}.
 \eeaa
Then $Y_{\t^*} = \wh X_{\t^*}$, $Y$ is an $\overline\cE^L$-supermartingale on $[0, \ch]$, and an $\overline\cE^L$-martingale on $[0, \t^*]$. Consequently, $\t^*$ is an optimal stopping time. 
\end{thm}

We define $\cP_\infty:=\cup_{L>0} \cP_L$. Our test processes will be smooth in the following sense.  
\begin {defn} (i) Let $u\in  C^0(\L,\dbK)$, its right time-derivative $\pa_t u$, if it exists, is defined as: 
\beaa
&&\pa_t u(t,\o):=\displaystyle \lim_{\d\downarrow 0}\frac{u(t+\d,\o_{. \wedge t})-u(t,\o)}{\d},\mbox{  for  } (t,\o)\in[0,T)\times \O,\\
&&\pa_t u(T,\o):=\displaystyle \lim_{t\uparrow T}\pa_t u(t,\o)\mbox{  for  } \o\in\O.
\eeaa
(ii) We say that $u$ is in  $C^{1,2}(\L)$ if $u,\pa_t u\in C^0(\L)$ and there exist $\pa_\o u \in C^0(\L,\dbR^d)$, $\pa_{\o\o} u \in C^0(\L,\dbS^d)$ such that the process $u_t:=u(t,B)$, verify :
$$d u_t=\pa_t u_t dt+ \frac12 \pa_{\o\o} u_t :d\la B\ra_t  +\pa_\o u_t  dB_t,\q \dbP\mbox{-a.s.},$$
for every $\dbP\in\cP_\infty$. 
\end{defn}
The previous notation $u_t:=u(t,B)$ will be our convention. If we compose a functional $u$ with the canonical process $B$, we will use the notation $u_t$. If we need to compose $u$ with other processes (for example $X$), we will explicitly write $u(t,X)$.
\begin{rem}
{\rm The requirement of continuity of the derivatives and the fact that the support of $\dbP_0$ is $\O$ show that if $u\in C^{1,2}({\L})$ then its derivatives are uniquely defined. \qed
}
\end{rem}
\subsection{Second order reflected BSDEs}
We refer to \cite{MPZ} for various properties of second order reflected BSDEs(2RBSDEs). The 2RBSDEs that we study have $3$ components :

- A final condition : $\xi: \O\rightarrow \dbR$. 

- A generator : $G:\L\times \dbR\times \dbR^d \times \dbS^d\rightarrow \dbR$.

- An obstacle : $h:\L\rightarrow \dbR$.

\subsubsection{The Generator}
Let $\cK$ be a measurable set with its sigma algebra $\cM_\cK$ and 2 mappings :
\beaa
&&F:\L\times\dbR\times\dbR^d\times \cK \rightarrow\dbR\\
&&\si:\L\times\cK \rightarrow\dbS^d.
\eeaa
We consider the following generator 
\bea\label{defn_generator}
&&G :\L\times\dbR\times\dbR^d\times \dbS^d \rightarrow \dbR\\
&&G(t,\o,y,z,\g)=\sup_{k\in\cK}\Big[\frac12 \sigma(t,\o,k)^2:\g+F(t,\o,y,\sigma(t,\o,k)z,k)\Big]
\eea

We make the following assumption on the data of the 2RBSDE:
\begin{assum}
\label{assumption_data} There exist $L_0,M_0\geq 0$ and $\rho_0$ a modulus of continuity with at most polynomial growth verifying the following points.

(i) \underline{Boundedness} : $\xi$, $h$, and $F(.,0,{\bf 0},.)$ are bounded by $M_0$. 

(ii) \underline{Assumptions on $F$ and $G$}: $F(.,y,z,k)$ and $G(.,y,z,\g)$ are right continuous under $d_\infty$ metric in the sense of the definition \reff{defn-rightcont}. $F(t,\o,.,.,k)$ is Lipschitz continuous in $(y,z)$ with Lipschitz constant $L_0$.

(iii)\underline{ Assumption on $h$}: $h$ is uniformly continuous under $d_\infty$ with modulus of continuity $\rho_0$. 

(iv)\underline{ Assumption on $\xi$} : $\xi$ is uniformly continuous under the $||.||_T$ norm with modulus of continuity $\rho_0$. For all $\o\in\O$, $\xi(\o)\geq h(T,\o)$.

(v)\underline{ Assumption on $\si$} : For all $(t,\o)\in \L$, $\inf_{k\in\cK}\si(t,\o,k)>0$, $|\sigma(t,\o,k)|\leq\sqrt{2L_0}$, $\si(t,.,k)$ is Lipschitz continuous with Lipschitz constant $L_0$, and $\si(.,k)$ is right continuous under $d_\infty$. 
\qed

\end{assum}

We will also need the following additional assumption for our wellposedness results. 
\begin{assum}
\label{assumption_data_add} $\si$ does not depend on $(t,\o)$, $F(.,y,z,k)$ is uniformly continuous with modulus of continuity $\rho_0$.
\qed
\end{assum}
\begin{rem}\label{rem_additional_assumption}
{\rm The Assumption \reff{assumption_data_add} will only be used to prove the Lemma \ref{equality} and under this additional assumption the operator $G$ is uniformly non-degenerate in $\g$.}
\end{rem}
\section{Introduction of the value functional}
\label{sect-intro-value}
For $t\in[0,T]$, we denote by $\cK^t$, the set of $\dbF^t$-progressively measurable and $\cK$ valued processes. Under the assumptions \reff{assumption_data}, for fixed $(t,\o)\in\L$, and $k\in \cK^t$, by the Lipschitz continuity of $\si$ in $\o$, there exists a unique strong solution $X^{t,\o,k}$ of the following equation under $\dbP^t_0$:
\bea\label{defn_X}
X^{t,\o,k}_s=\int_t^s \si^{t,\o}(r,X^{t,\o,k},k_r)dB^t_r,\mbox{ for } s\in[t,T].
\eea

Additionally, by the classical estimates on SDEs, for $(t,\o),(t,\o')\in \L$ :
\bea\label{estimate_sde}
&&\dbE^t_0 \left[(||X^{t,\o,k}||_T^t)^p\right]\leq C_p,\mbox{ for all }p>0,\notag\\
&&\dbE^t_0 \left[(||X^{t,\o,k}-X^{t,\o',k}||_T^t)^2\right]\leq C||\o-\o'||_t^2.
\eea

At the previous inequality, as it will be the case in the sequel, $C$ is a constant that may change from line to line, however it only depends on $d,M_0,T,L_0$, and $ \rho_0$.

We define $\dbP^{t,\o,k}:=\dbP^t_0 \circ (X^{t,\o,k})^{-1}\in\cP^t_{L_0}$. Notice that, the lemma 2.2 of \cite{STZ_target} shows that there exits a mapping $\tilde k\in\dbL^0(\L;\cK)$ such that $\tilde k(s,X^{t,\o,k})=k(s,B^t)$, $ds\times \dbP^t_0$-a.s. By rewriting \reff{defn_X} under $\dbP^t_0$:
\bea\label{defn_X2}
X^{t,\o,k}_s=\int_t^s \si^{t,\o}(r,X^{t,\o,k},\tilde k(r,X^{t,\o,k}))dB^t_r,\mbox{ for } s\in[t,T].
\eea

Therefore, $\{\sigma^{t,\o} _r(\tilde k_r)^{-1} d B_r^{t}\}_{r\in[t,T]}$(recall that  $\sigma^{t,\o}_r(\tilde k_r)=\sigma^{t,\o} (r,B^t,\tilde k(r,B^t))$ ) is the increment of a Brownian motion under $\dbP^{t,\o,k}$. Hence, for fixed $\t\in\cT^t$ and $\dbF^t_\t$ measurable and bounded random variable $\zeta$, one can define $(\cY^{t,\o,k}_s(\t,\zeta),\cZ^{t,\o,k}_s(\t,\zeta), \cK^{t,\o,k}_s(\t,\zeta))_{s \in [t,\t]}$  solution to the reflected BSDE on $\O^t$, with data $(F^{t,\o}_s(.,.,\tilde k_s),h^{t,\o},\zeta)$ under $\dbP^{t,\o,k}$:
\bea\label{RBSDE-general}
&&\cY_s^{t,\o,k}=\zeta(B^{t})+ \int_s^\t {F_r^{t,\o} ( \cY_r^{t,\o,k},\cZ_r^{t,\o,k},\tilde k_r) dr }\\
&&\q\q\q-\int_s^\t {(\cZ_r^{t,\o,k})^*\sigma_r^{t,\o} (\tilde k_r)^{-1} d B_r^{t}}+\cK_\t^{t,\o,k}-\cK_s^{t,\o,k},\notag\\
&& \cY_s^{t,\o,k} \geq h^{t,\o}_s, \mbox{ for all }s\in [t,\t],\notag\\
&&(\cK_s^{t,\o,k})_{s\in [t,\t]} \mbox { is increasing in s, } \cK_t^{t,\o,k}=0 \mbox{ and }\left[\cY_s^{t,\o,k}-h^{t,\o}_s\right] d\cK_s^{t,\o,k}=0.\notag
\eea
When $(\t,\zeta)=(T,\xi)$, we denote 
\bea\label{RBSDE}
(Y^{t,\o,k}_s,Z^{t,\o,k}_s, K^{t,\o,k}_s)=(\cY^{t,\o,k}_s(T,\xi),\cZ^{t,\o,k}_s(T,\xi), \cK^{t,\o,k}_s(T,\xi)).
\eea
To make easier our notations, we also define the following reflected BSDE, under $\dbP^t_0$:
\bea\label{RBSDE_tilde}
&&\tilde Y_s^{t,\o,k}=\xi^{t,\o}(X^{t,\o,k})+ \int_s^T {F^{t,\o} (r, X^{t,\o,k}, \tilde Y_r^{t,\o,k},\tilde Z_r^{t,\o,k},\tilde k(r,X^{t,\o,k})) dr }\\
&&\q\q\q-\int_s^T {(\tilde Z_r^{t,\o,k})^* d B_r^{t}}+\tilde K_T^{t,\o,k}-\tilde K_s^{t,\o,k},\notag\\
&& \tilde Y_s^{t,\o,k} \geq h^{t,\o} (s,X^{t,\o,k}), \mbox{ for all }s\in [t,T],\notag\\
&&(\tilde K_s^{t,\o,k})_{s\in [t,T]} \mbox { is increasing in s, } \tilde K_t^{t,\o,k}=0 \mbox{ and }\left[\tilde Y_s^{t,\o,k}-h^{t,\o} (s,X^{t,\o,k})\right] d\tilde K_s^{t,\o,k}=0.\notag
\eea
For all $s\in [t,T]$, $Y_s^{t,\o,k}$ and $\tilde Y_s^{t,\o,k}$ are $\overline {\dbF^t}^{\dbP^{t,k}}_s$ measurable, so $Y_t^{t,\o,k}$ and $\tilde Y_t^{t,\o,k}$ are constant. Additionally, the family: $$(\xi^{t,\o}(B^t),F^{t,\o}_s(y,z,\tilde k_s),h^{t,\o}_s,\sigma^{t,\o}_s (\tilde k _s)^{-1} d B_s^{t})_{s\in[t,T]}$$ under $\dbP^{t,\o,k}$  has the same distribution as the family  $$(\xi^{t,\o}(X^{t,\o,k}),F^{t,\o}(s,X^{t,\o,k},y,z,\tilde k(s,X^{t,\o,k})),h^{t,\o}(s,X^{t,\o,k}),d B_s^{t})_{s\in[t,T]}$$ under $\dbP_0^{t}$. So  $Y_t^{t,\o,k}=\tilde Y_t^{t,\o,k}$.  We define the following process which is our value functional of interest :
\bea\label{value_function}
u^0 (t,\o):= \sup_{k\in\cK^t} Y_t^{t,\o,k}=\sup_{k\in\cK^t} \tilde Y_t^{t,\o,k},\mbox{ for } (t,\o)\in \L.
\eea

\subsection{Regularity of the value functional}
\begin{prop}\label{regularity_value}
$u^0$ is bounded and uniformly continuous under the $d_\infty$ metric in $\L$. 
\end{prop}
\proof Under assumptions (\ref{assumption_data}), the data of the problem verifies the assumptions of \cite{EKKPQ} which gives the following a priori estimate :
\beaa
\dbE^{t,\o,k}\left (\sup_{t\leq s \leq T} |Y_s^{t,\o,k}|^2 +\int_t^T |Z_r^{t,\o,k}|^2 dr+ (K_T^{t,\o,k})^2 \right) \leq C.
\eeaa
Additionally, for $(t,\o),(t,\o')\in\L$, notice that $||\o\otimes_t X^{t,\o,k}-\o'\otimes_t X^{t,\o'.k}||_T\leq ||\o-\o'||_t+||X^{t,\o,k}- X^{t,\o',k}||^t_T$, therefore under our boundedness and regularity assumptions the estimates in \cite{EKKPQ} gives : 
\bea\label{ekkpq_difference_estimate}
&&|\tilde Y_t^{t,\o,k}-\tilde Y_t^{t,\o',k}|^2 \notag\\
&&\leq C \dbE_0^{t}\left[ |\xi^{t,\o}(X^{t,\o,k})-\xi^{t,\o'}(X^{t,\o',k})|^2\right]\notag\\
&&+C\dbE^t_0 \left[\int_t^T {|F^{t,\o} (s,X^{t,\o,k}, Y_s^{t,\o,k},Z_s^{t,\o,k}, k_s)-F^{t,\o'} (s,X^{t,\o',k}, Y_s^{t,\o,k},Z_s^{t,\o,k},k_s)|^2 ds}\right]\notag\\
&&+C\dbE^{t}_0 \left [\sup_{t\leq s\leq T}| h^{t,\o}(s,X^{t,\o,k})-h^{t,\o'}(s,X^{t,\o',k})|^2\right]^{1/2}\notag\\
&&\leq C \dbE_0^{t}\left[ \rho_0^2(||\o'-\o||_{t}+||X^{t,\o,k}-X^{t,\o',k}||^t_T)\right].
\eea
$\rho_0$ has at most polynomial growth, denote $p_0>0$ this growth power. For fixed $\d>0$, we can estimates the difference $ |\tilde Y_t^{t,\o,k}-\tilde Y_t^{t,\o',k}|$  as follows :
\bea\label{var_omega_uk}
&&|\tilde Y_t^{t,\o,k}-\tilde Y_t^{t,\o',k}|^2\\
&&\leq C\dbE^{t}_0 \left(\rho_0^2(||\o'-\o||_{t}+||X^{t,\o,k}-X^{t,\o',k}||^t_T)\1_{\{||X^{t,\o,k}-X^{t,\o',k}||^t_T> \d\}}\right)\notag\\
&&+C\dbE^{t}_0 \left(\rho_0^2(||\o'-\o||_{t}+\d)\1_{\{||X^{t,\o,k}-X^{t,\o',k}||^t_T\leq \d\}}\right)\notag\\
&&\leq C \sqrt{\frac{\dbE^{t}_0 \left(||X^{t,\o,k}-X^{t,\o',k}||^t_T\right)}{\d} }\left[1+ ||\o-\o'||^{2p_0}_t \right]+  C\rho_0^2(||\o'-\o||_{t}+\d)\notag\\
&&\leq C\left(\sqrt{\frac{||\o'-\o||_{t}}{\d}}\left[1+ ||\o-\o'||^{2p_0}_t \right]+\rho_0^2(||\o'-\o||_{t}+\d)\right).\notag
\eea
If we choose $\d:=\sqrt{||\o-\o'||_t}$, then the last line becomes a modulus of continuity $\rho_1$ with at most polynomial growth.

First of all, the previous estimates gives that $Y_t^{t,\o,k}$ is bounded by a constant that only depends on $M_0,T,L_0$, and $ \rho_0$. With a passage to supremum in $k$, we see that $u^0$ is bounded. Additionally  
\bea\label{regularity_o}
|u^0(t,\o)-u^0(t,\o')|\leq \sup_{k\in\cK^t}|\tilde Y_t^{t,\o,k}-\tilde Y_t^{t,\o',k}|\leq \rho_1(||\o-\o'||_t),
\eea
which show that for fixed $t$, $u^0$ is uniformly continuous in $\o$ uniformly in $t$.

Fix $0\leq t\leq t_1\leq T$. Given the uniform continuity of $u^0$ in $\o$ for fixed times, one can proceed as in Lemma 4.1 of \cite{ETZ0} to obtain the following dynamic programming principle at deterministic times :
\bea\label{DPP-deterministic}
u^0(t,\o)= \sup_{k\in \cK^t} \cY^{t,\o,k}_t(t_1,u^0(t_1,\o\otimes_t B^t)),
\eea
where $\cY^{t,\o,k}(t_1,u^0(t_1,\o\otimes_t B^t))$ (denoted only by $\cY^{t,\o,k}$ for simplicity at this section) is defined at \reff{RBSDE-general}.
We estimate the variation in time for $k\in\cK^t$ and under $\dbP^{t,\o,k}$ 
\beaa
u^0(t_1,\o)-u^0(t,\o)&=&u^0(t_1,\o)-u^0(t_1,\o\otimes_t B^t)-\int_t^{t_1}F^{t,\o}(r,B^t,\cY^{t,\o,k}_r,\cZ^{t,\o,k}_r,\tilde k_r)dr\\
&+&\int_t^{t_1}\cZ^{t,\o,k}_r \cdot (\si^{t,\o}_r(\tilde k_r))^{-1}dB^t_r-\cK^{t,\o,k}_{t_1}+\cY^{t,\o,k}_t-u^0(t,\o)
\eeaa
We take the expectation under $\dbP^{t,\o,k}$ to have : 
\beaa
u^0(t_1,\o)-u^0(t,\o)&\leq&\dbE^{\dbP{t,\o,k}}\left[\rho_1 (||\o-\o\otimes_t B^t||_{t_1})\right]+|\cY^{t,\o,k}_t-u^0(t,\o)|\\
&+&C(t_1-t)+ L_0 \int_t^{t_1}\dbE^{\dbP^{t,\o,k}}\left[|\cZ^{t,\o,k}_r|\right]dr.
\eeaa
Finally by taking a sequence $k_n$ such that $\cY^{t,\o,k_n}_t\rightarrow u^0(t,\o)$ and  using  estimates on the RBSDEs we have:
\beaa
u^0(t_1,\o)-u^0(t,\o)&\leq& C\sqrt{t_1-t}+\rho(||\o_{.\wedge t}-\o||_{t_1}),
\eeaa
for some modulus of continuity $\rho$. 

We define the following optimal stopping time for $\cY^{t,\o,k}$,
$D^{t,\o,k}:=\inf\{s\in[ t,t_1] : \cY^{t,\o,k}_s=h^{t,\o,k}_s\}\wedge t_1. $
Then 
\beaa
u^0(t_1,\o)-\cY^{t,\o,k}_t&=&\dbE^{\dbP^{t,\o,k}}\left[u^0(t_1,\o)-(u^0)^{t,\o}_{t_1}+\1_{\{D^{t,\o,k}<t_1\}}((u^0)^{t,\o}_{t_1}-h^{t,\o}_{D^{t,\o,k}})\right.\\
&-&\left. \int_s^{t_1\wedge D^{t,\o,k}} F^{t,\o}(r, B^t_\cd,  \cY_r,\cZ_r, \tilde k_r) dr\right]\\
&\geq&\dbE^{\dbP^{t,\o,k}}\left[-\rho_1(||\o-\o\otimes_t B^t||_{t_1})+\1_{\{D^{t,\o,k}<t_1\}}(h^{t,\o}_{t_1}-h^{t,\o}_{D^{t,\o,k}})\right.\\
&-&\left. \int_s^{t_1} |F^{t,\o}(r, B^t_\cd,  \cY_r,\cZ_r, \tilde k_r)| dr\right]
\eeaa
Recall that $h$ is uniformly continuous hence we can bound the term $(h^{t,\o}_{t_1}-h^{t,\o}_{D^{t,\o,k}})$ on the event $\{D^{t,\o,k}<t_1\}$. Therefore, we can control the right hand side uniformly in $k$. Finally combining all the previous results, we obtain that there is a modulus of continuity $\tilde \rho_0$, which only depends on, $M_0,L_0,\rho_0,T $, such that 
\bea\label{unif-cont-u0}
|u^0(t,\o)-u^0(t',\o')|\leq \tilde \rho_0 (d_\infty((t,\o),(t',\o'))).
\eea
\qed

\section{Viscosity solutions to path-dependent PDEs}

For any $L\ge 0$ and $(t,\o)\in [0,T)\times \O$, and $u\in \underline \cU$, define:
\bea
\label{cA}
\left.\ba{lll}
\underline\cA^{L} u (t,\o)
:=
\Big\{\f \in C^{1,2}_b(\L^t):~\mbox{there exists $\ch\in\cH^t$ such that}
\\
\hspace{40mm}
0=\f(t,{\bf 0})-u(t,\o)
=\underline\cS^{L}_t\big[(\f- u^{t,\o})_{.\wedge\ch}
\big](\bf 0)
\Big\};\\
\ea\right.
\eea
and for all, $u\in \overline \cU$:
\bea
\left.\ba{lll}
\\ \overline\cA^{L}u (t,\o)
:=
\Big\{\f \in C^{1,2}_b(\L^t): ~\mbox{there exists $\ch\in\cH^t$ such that}
\\
\hspace{40mm}
0=\f(t,{\bf 0})-u(t,\o)
=\overline\cS^{L}_t\big[(\f- u^{t,\o})_{.\wedge\ch}
\big](\bf 0)
\Big\}.
\ea\right.
\eea
These sets are the equivalents of sub/superjets in our theory. 

\subsection{The PPDE}
For fixed $(t,\tilde \o)\in \L$, we define the differential operator $\cL^{t,\tilde \o}$ on $C^{1,2}(\L^t)$: 
\beaa
&&\mbox{for }\phi \in C^{1,2}(\L^t)\mbox{ and }(s,\o)\in \L^t\notag\\
&&\cL^{t,\tilde \o} \phi(s, \o) :=-\pa_t \phi (s,\o)-G^{t,\tilde\o}(s,\o,\phi (s,\o), \pa_\o\phi(s,\o),\pa_{\o\o}\f(s,\o)).
\eeaa
When $t=0$ the operator is simply written $\cL$.
The functional $u^0$ defined by \reff{value_function} is related, as it is the case in the Markovian case, to the following PPDE :
\bea
\label{PPDE}
&&\min\{\cL u (t,\o);(u-h)(t,\o)\}=0,\mbox{  for all }(t,\o)\in [0,T)\times \O,\\
&&u(T,\o)=\xi(\o) \mbox{  for all }\o\in\O.
\eea 
\subsection{Viscosity solution of PPDEs}
We give the following definition of viscosity solution. 
\begin{defn}
\label{defn-viscosity}
(i) For any $L\ge 0$, we say $u\in \overline \cU$ is a viscosity $L$-supersolution of PPDE \reff{PPDE} if, for any $(t,\o)\in [0, T)\times \O$ and any $\f \in \overline\cA^{L}u  (t,\o)$, it holds that
\beaa
&&u(t,\o)-h(t,\o) \ge  0,\mbox{ and } (\cL^{t,\o}\f)(t,{\bf 0}) \ge  0,\\
&&\mbox{or equivalently } \min\{\cL^{t,\o} \phi(t,{\bf 0}); u(t,\o)-h(t,\o) \}\geq 0.
\eeaa

(ii) We say $u\in\underline \cU$ is a viscosity $L$-subsolution of PPDE \reff{PPDE} if, for any $(t,\o)\in [0, T)\times \O$ such that $u(t,\o)-h(t,\o)>0$ and any $\f \in \underline\cA^{L} u (t,\o)$, it holds that
\beaa
&&(\cL^{t,\o}\f)(t,{\bf 0}) \le  0.
\eeaa

(iii) We say $u$ is a viscosity subsolution (resp. supersolution) of PPDE \reff{PPDE} if $u$ is viscosity $L$-subsolution (resp. $L$-supersolution) of PPDE \reff{PPDE} for some $L\ge 0$.

(iv) We say $u$ is a viscosity solution of PPDE \reff{PPDE} if it is both a viscosity subsolution and a viscosity supersolution.
\end{defn}
\begin{rem}\label{increasing_prop}
{\rm For $0\leq L_1\leq L_2$ and $(t,\o)\in [0,T)\times \O$, we have $\underline\cE^{L_2}_t[.]\leq\underline \cE^{L_1}_t[.]$ and $\underline \cA^{L_2} u(t,\o)\subset \underline \cA^{L_1} u(t,\o)$. If $u$ is a viscosity $L_1$-subsolution then $u$ is a viscosity $L_2$-subsolution. Same statement also holds for supersolutions.\qed}
\end{rem}
\begin{rem}
\label{rem-cA}
{\rm The definition of viscosity solution property is local in the following sense. For any $(t,\o) \in [0,T)\times \O$, to check the viscosity property of $u$ at $(t,\o)$, it suffices to know the value of $u^{t,\o}$ on $[t, \ch^t_\d]$ for an arbitrarily small $\d>0$. The hitting times $\ch^t_\d$ are our tools of localization. \qed
}
\end{rem}

\begin{rem}
\label{visco_sol_strict}
{\rm We have some flexibility to choose the set of test functionals. All the results in this paper still hold true if we replace the $\underline{\cA}^Lu$  with the ${\underline\cA'}^{L}\!u$ 
\bea
\label{cA1}
{\underline\cA'}^{L}\!u (t,\o)
&:=&
\Big\{\f\!\in\!C^{1,2}(\L^{\!t}):
      \exists\;\ch\in\cH^t~\mbox{such that, for all}~\t'\in\cT^t_+,\nonumber\\
&&\qq   (\f-u^{t,\o})_t({\bf 0}) \;=\; 0
       \;<\;\underline\cE^{L}_t
        \!\big[(\f- u^{t,\o})_{\t'\wedge\ch}
        \big]
\Big\}.
\eea
\qed}
\end{rem}

\subsection{Consistency}
\begin{prop}
Assume $u\in C^{1,2}(\L,\dbR)$ then $u$ is a viscosity subsolution(respectively, supersolution) of the PPDE (\ref{PPDE}) if and only if $u$ is a classical subsolution(respectively, supersolution) of the same equation. 
\end{prop}
\proof 
We only prove the subsolution case. A similar proof also holds for supersolutions. Assume that $u$ is a viscosity $L$-subsolution and take $(t,\o)\in[0,T)\times \O$ and $u(t,\o)-h(t,\o)>0$. Choosing $\phi=:u$ and $\ch:=T\in\cH^t$, clearly $\phi\in \underline \cA^L u(t,\o)$ so $\cL^{t,\o} \phi (t,\o)=\cL u(t,\o)\leq 0$.

For the reverse implication, assume that $u$ is a classical subsolution and it is not a viscosity subsolution, a fortiori it is not a viscosity $L_0$-subsolution. Therefore, there exist $(t,\o)\in[0,T)\times \O$ and $\phi\in\underline\cA^{L_0} u (t,\o)$ with the associated $\ch\in\cH^t$ such that $(u-h)(t,\o)>0$ and  $c:=\cL^{t,\o} \phi(t,{\bf 0 })>0$. The processes $G(.,\phi_t,\pa_\o\phi_t,\pa_{\o\o}\f_t),\phi_.,u^{t,\o}_.$ are right continuous under the $d_\infty$ metric, so there exist $\d>0$ such that for $(s,\tilde\o)\in[t,\ch^t_\d]\times \O^t$, the following  inequalities holds:
\bea\label{construction_tilde_tau1}
&&|G(t,\o,\phi_t,\pa_\o\phi_t,\pa_{\o\o}\phi_t)-G(s,\o\otimes_t\tilde \o,\phi_t,\pa_\o\phi_t,\pa_{\o\o}\phi_t)|\leq c/4,\\
&&|\pa_t\phi_t-\pa_t\phi_s|+L_0 |\phi_t-\phi_s|+M_0L_0|\pa_\o \phi_t-\pa_\o\phi_s|+L_0|\phi_s-u^{t,\o}_s|+L_0|\pa_{\o\o}\phi_t-\pa_{\o\o}\phi_s|\leq c/4\notag
\eea
Then $\cL^{t,\o} \phi_s-L_0 |\phi_s-u_s^{t,\o}|\geq c/2$  for $s\in[t,\ch^t_\d]$. $u$ is a subsolution of the PPDE \reff{PPDE} and the data is right continuous under $d_\infty$, so we can choose a constant process $k\in\cK^t$ and $\d>0$ small enough such that for all $s\in[t,\ch^t_\d]$ :
\bea\label{control_k}
&&\pa_t u_s^{t,\o} +\frac12 \pa_{\o\o} u_s^{t,\o} :(\si_s(k_s))^2+F^{t,\o}_s(u^{t,\o}_s,\si_s(k_s) \pa_\o u^{t,\o}_s ,k_s)\geq-c/4.
\eea
Notice that for $k$ constant the equation~\reff{defn_X} has strong solutions. Applying It\^o's formula under $\dbP^{t,\o,k}$:
\beaa
&&0= (\phi -u^{t,\o})_t=(\phi -u^{t,\o} )_{\ch^t_\d}\\
&& -\int_t^{{\ch^t_\d}} \pa_t(\phi-u^{t,\o})_s+\frac{1}{2}\pa_{\o\o} (\phi -u^{t,\o})_s:\si_s(k_s)^2ds -\int_t^{{\ch^t_\d}} \pa_\o(\phi- u^{t,\o})^*_s  dB_s^t\\
&&\geq (\phi -u^{t,\o} )_{{\ch^t_\d}}+\int_t^{{\ch^t_\d}} {\cL^{t,\o}\phi_s ds}+\int_t^{{\ch^t_\d}} {-c/4 ds}\\
&&+\int_t^{{\ch^t_\d}} r_s(\phi-u^{t,\o})_s +\a_s^* \si_s(k_s) \pa_\o(\phi-u^{t,\o})_s ds -\int_t^{{\ch^t_\d}}\pa_\o(\phi- u^{t,\o})^*_s  dB_s^t\\
\eeaa
For some $r_s$ and $\a_s$  progressively measurable $r\in\dbR$, $\a\in\dbR^d$, with $|r_s|\leq L_0$ and $|a_s|\leq L_0$ . Therefore  :
\beaa
&&0\geq (\phi -u^{t,\o} )_{{\ch^t_d}}+\int_t^{{\ch^t_d}} {\cL^{t,\o}\phi_s +r_s(\phi-u^{t,\o})_s ds}\\&&-\frac{(\ch^t_d-t) c}{4 } -\int_t^{{\ch^t_d}}\pa_\o(\phi- u^{t,\o})^*_s  (dB_s^t-\si_s(k_s) \a_s ds)\\
&&\geq  (\phi -u )_{{\ch^t_d}}+\frac{(\ch^t_d-t) c}{4 } -\int_t^{{\ch^t_d}}\pa_\o(\phi- u^{t,\o})^*_s  (dB_s^t-\si_s(k_s) \a_s ds)
\eeaa
Notice that by Girsanov's theorem, there exists $\dbP\in \cP^t_{L_0}$ equivalent to $\dbP^{t,\o,k}$ such that  the last  integral is a martingale under $\dbP$. Therefore, we have the following inequalities that contradicts the assumption that $\phi\in\underline \cA^{L_0}u (t,\o)$ : 
\beaa
&&0>-\frac{c}{4} \dbE^{\dbP} [\ch^t_\d-t]\geq\dbE^{\dbP} [(\phi -u )_{\ch^t_\d}]\geq \underline \cE^{L_0}_t[(\phi -u )_{\ch^t_\d}].
\eeaa
\qed 

\subsection{A change of variable formula}
We will need the following change of variable formula in our subsequent analysis. 
\begin{prop}
\label{change_var}
Let $C,\l,\mu\in\dbR$ be constants and $u\in\underline \cU$ then $u$ is a viscosity $L$-subsolution of the PPDE \reff{PPDE} with data $(G,\xi, h)$ if and only if $ u'_t:=e^{\l t} u_t+C e^{\mu t}t$ is a viscosity $L$-subsolution of the PPDE \reff{PPDE} with data $( G',\xi',h')$ where :
\beaa
&& G'(t,\o,y,z,\g):=-Ce^{\mu t} (1+(\mu-\l)t)-\l y + e^{\l t} G(t,\o,e^{-\l t}y-C e^{(\mu-\l)t}t,e^{-\l t}z,e^{-\l t}\g ),\\
&& \xi':=e^{ \l T} \xi+C e^{\mu T}T,\\
&&h'_t:=e^{\l t}h_t +Ce^{\mu t} t.
\eeaa
The same statement holds also for $L$-supersolutions. 
\end{prop}
\proof
We will only prove the subsolution case. Assume that $u$ is a viscosity $L$-subsolution with data $(G,\xi, h)$. We want to show that $ u'$ is a viscosity L-subsolution with data $ (G',\xi', h')$. Take $(t,\o)\in [0,T)\times \O$  such that $u'(t,\o)- h'(t,\o)>0$ (notice that this is equivalent to $u(t,\o)-h(t,\o)>0$) and $\phi'\in\underline \cA^L  u' (t,\o)$, with the corresponding hitting time $\ch\in\cH^t$.\\
For fixed $\e>0$, we define $$\phi_s^\e := e^{-\l s}  \phi'_s -Ce^{(\mu-\l)s} s+\e (s-t).$$
Notice that $\phi'_t= u' _t =e^{\l t} u_t +Ce^{\mu t}t$ and for $s\geq t$:
\beaa
&&\phi_s^\e -u_s^{t,\o} -e^{-\l t} ( \phi'_s- (u'_s)^{t,\o}) \\
&&= (e^{-\l s}-e^{-\l t })(( \phi'_s-C e^{\mu s} s) -( \phi'_t-C e^{\mu t} t))+ (e^{\l (s-t)}-1)(u_s^{t,\o}-u_t^{t,\o})\\ 
&&+ (e^{\l (t-s)}+ e^{\l (s-t)}-2) u_t^{t,\o}+ \e (s-t)
\eeaa
There exists a constant $K>1$ which may depend on $\l,t,T$  but not in $s\in(t,T]$, and $\e>0$ such that :
\beaa
&&0\leq |(e^{-\l s}-e^{-\l t })|\leq K(s-t)\\
&&0\leq e^{\l (s-t)}-1\leq K(s-t)\\
&&0\leq |e^{\l (t-s)}+ e^{\l (s-t)}-2|\leq K(s-t)^2
\eeaa 
Additionally $u\in\underline \cU$ and $\phi$ is continuous in under the $d_\infty$ metric, so there exist $\d$(depending on $\e$) such that on $[t,\ch^t_\d]$ :
\beaa
&&R_0:=1\vee \sup_s u_s<\infty,\\
&&u_s^{t,\o}-u_t^{t,\o}\geq\frac{-\e}{3K R_0}\\
&&|( \phi'_s-C e^{\mu s} s) -( \phi'_t-C e^{\mu t} t)|\leq \frac{\e}{3K R_0}\\
&&0\leq s-t\leq \frac{\e}{3K R_0}.
\eeaa
Combining the previous inequalities :
\beaa
\phi_s^\e -u_s^{t,\o} -e^{-\l t} ( \phi'_s- (u'_s)^{t,\o})\geq -\e(s-t)+\e(s-t)\geq 0.
\eeaa

Then for all $\t\in\cT^t$, such that $\t\leq\ch^t_\d$ it holds that :
$$\phi_\t^\e-u^{t,\o}_\t \geq e^{-\l t }( \phi'_\t - (u')^{t,\o}_\t),$$
therefore :
$$\underline \cE^{L}_t [\phi_\t^\e -u^{t,\o}_\t]\geq e^{-\l t}\underline \cE^{L}_t [ \phi'_\t- (u')^{t,\o}_\t]\geq 0=\phi_t^\e-u_t^{t,\o}.$$
Which shows that $\phi^\e \in \underline \cA   u(t,\o)$, and $u(t,\o)-h(t,\o)>0$, by the definition of viscosity subsolutions,
 $0\geq \cL^{t,\o} \phi^\e $. Taking the limit as $\e$ goes to $0$, we obtain :
\beaa
&&0\geq \cL^{t,\o} \phi =-\pa_t \phi'_t -G' (t,\o,\phi'_t,\pa_\o\phi'_t,\pa_{\o\o}\phi'_t).
\eeaa
\qed

\begin{rem}
\label{rem_change}

{\rm Notice that after the change of variable the function $G'$ can be written as :
\beaa
&&G'(t,\o,y,z,\g)=\sup_{k\in\cK}\left\{\frac{\si(t,\o,k)^2:\g}{2}+F'(t,\o,y,\si(t,\o,k)z,k)\right\}\\
&&\mbox{where }F'(t,\o,y,z,k):=e^{\l t}F(t,\o,e^{-\l t}y-Ce^{(\mu-\l)t}t,e^{-\l t}z,k)-Ce^{\mu t}(1+(\mu-\l)t)-\l y
\eeaa
We will make the following choices for the constants: $\l=L_0+1$, $\mu=0$ and $C=-2 e^{(L_0+1) T}(\l+1)(M_0+1) $. With this change of variable, the data of the problem verify the following properties: 
\bea\label{monotonicity_f}
&&G'(t,\o,y,z,\g)\geq G'(t,\o,y+\d,z,\g)+\d, \mbox { for all } \d>0, and \mbox { any } (t,\o,y,z,\g),\\
&&F'(t,\o,y,z,k)\geq F'(t,\o,y+\d,z,k)+\d, \mbox { for all } \d>0, and \mbox { any } (t,\o,y,z,k),\notag\\
&& F'(t,\o, h'(t,\o),{\bf 0},k)=-C-\l e^{\l t}h_t +e^{\l t}F(t,\o,h(t,\o),{\bf 0},k)\geq 0\mbox { for all  } (t,\o)\in \L .\notag
\eea 
When needed, we will assume, that $F,G$ and $h$ verify (\ref{monotonicity_f}). This change of variable formula will be useful at subsection \reff{study_K}. }
\end{rem}

\subsection{Viscosity solution property of the value functional}
Before starting to study the viscosity solution property of $u^0$, we give the following dynamic programming principal on random times. Its proof is similar to the proof of Theorem 4.3 of \cite{ETZ0}.
With the notation introduced at \reff{RBSDE-general}, for all $\t\in\cT^t$, the following dynamic programming at stopping times holds
\bea\label{second_dpp}
u^0(t,\o)=\sup_{k\in\cK^t}\cY^{t,\o,k}_t(\t,u^0 (\t,\o\otimes_t B^t))
\eea

\begin{thm}
Under the assumptions \reff{assumption_data} on the data of the problem, the value functional $u^0$ defined at \reff{value_function} is viscosity solution of the PPDE \reff{PPDE}.
\end{thm}

\subsubsection{Subsolution property of the value functional}
We assume without loss of generality that $$G \mbox { and }F \mbox{ are increasing in }y.$$
We reason by contradiction by assuming that $u^0$ is not a viscosity $L_0$-subsolution, so there exist $(t,\o)\in[0,T)\times \O$ such that $u^0(t,\o)> h(t,\o)$, and $\phi\in\underline \cA^{L_0} u^0(t,\o)$, with the associated $\ch\in\cH^t$ and verifying: $$c=\min\{\cL^{t,\o} \phi (t,0),u^0(t,\o)-h(t,\o)\}>0.$$

Without loss of generality we will assume that $(t,\o)=(0,{\bf 0})$. Recall that $u^0$ and $h$ are uniformly continuous. Therefore there exist $\d>0$ such that for all $s\in[0,\ch_\d]$, $$|u^0_0-u^0_s|\leq c/4,\q|h_0-h_s|\leq c/4.$$ 
Denote, for $s\in [0,T]$ and $k\in\cK^0$ :
\beaa
&&  \d Y_s^k:=\phi_s - Y^{0,{\bf 0},k}_s,\q \d Z_s^k :=  \pa_\o\phi_s-(\si_s(k_s))^{-1} Z^{0,{\bf 0},k}_s\\
&& G_s(u^0,\pa_\o \phi):=G_s(u^0_s, \pa_\o\phi_s,\pa_{\o\o}\phi(s,B)),\q F^k_s(y,z):=F_s(y, z,\tilde k_s),\notag
\eeaa
By the continuity of $\phi$ and the right continuity of $G$, we can take $\d>0$ small enough, to have for $s\in[0,\ch_\d]$: $$-\pa_t\phi_s-G_s(u^0_s,\pa_\o \phi_s,\pa_{\o\o}\phi_s)\geq c/2.$$
$G$ is defined as the supremum in \reff{defn_generator}, then for all $k\in\cK^0$, the following inequality holds : 
$$-\pa_t\phi_s-\frac12 \pa_{\o\o}\phi_s:\si_s({k}_s)^2-F_s(u^0_s,\si_s({k}_s)\pa_\o \phi_s,\tilde k_s)\geq c/2.$$
For $k\in\cK^0$, we apply functional It\^o's formula to $\phi$ and use the definition of $Y^{0,{\bf 0},k}$ in \reff{RBSDE} to obtain under $\dbP^{0,{\bf 0},k}$ :
\beaa
d (\d Y^k)_s &=& [ \pa_t \phi_s +\frac12 \pa_{\o\o}\phi_s :\si_s(k_s)^2 +F_s^k(u^0,\pa_\o\phi_s)]ds\\
&&+[F_s^k(Y^{0,{\bf 0},k}_s,Z^{0,{\bf 0},k}_s)- F_s^k(u^0_s, \pa_\o\phi_s) ]ds+(\d Z^k_s)^*  dB_s+dK^{0,{\bf 0},k}_s.
\eeaa
Therefore   for all $k$, $\dbP^{0,{\bf 0},k}$-a.s. : 
\beaa
\d Y^k_{\ch_\d}-\d Y^k_{0} &=& \int_t^{\ch_d} [\pa_t \phi_s + \frac12 \pa_{\o\o}\phi(s,B^t) :\si_s(k_s)^2 +F_s^k(u^0,\phi)]ds\\
&&+\int_0^{\ch_d}(F_s^k(Y^{0,{\bf 0},k}_s,Z^{0,{\bf 0},k}_s)- F_s^k(u^0_s, \phi))ds +\int_0^{\ch_d}(\d Z^k_s)^*  dB_s+K^{0,{\bf 0},k}_{\ch_d}\\
&&\leq \frac{-c \ch_\d}{2}+\int_t^{\ch_\d}(F_s^k(Y^{0,{\bf 0},k}_s,Z^{0,{\bf 0},k}_s)- F_s^k(u^0_s, \phi))ds +\int_t^{\ch_\d}(\d Z^k_s)^*  dB_s+K^{0,{\bf 0},k}_{\ch_\d}
\eeaa
We have assumed that $F$ is increasing in $y$ and $u^0_s\geq Y_s^{0,{\bf 0},k}$ therefore for all $k$, $\dbP^{0,{\bf 0},k}$-a.s.:
\beaa
&&(\phi-u^0)_{\ch_\d}-(u^0-Y^{0,{\bf 0},k})_{0}=(\phi-u^0)_{\ch_\d}-(\phi-Y^{0,{\bf 0},k})_{0}\\
&&\leq (\phi-Y^{0,{\bf 0},k})_{\ch_\d}-(\phi-Y^{0,{\bf 0},k})_{0}=\d Y^k_{\ch_\d}-\d Y^k_{0}\\
&&\leq\frac{-c \ch_\d}{2}+\int_0^{\ch_\d}(F_s(u^0_s, Z_s^{0,{\bf 0},k})- F_s^k(u^0_s, \phi))ds +\int_0^{\ch_\d}(\d Z^k_s)^*  dB_s+K^{0,{\bf 0},k}_{\ch_\d}\\
&&=\frac{-c \ch_\d}{2}+\int_0^{\ch_\d}(\d Z^k_s)^* ( dB_s+ \si_s(k_s)\a_s ds)+K^{t,\o,k}_{\ch^t_\d}
\eeaa
where $|\a_s|\le L_0$. 

By the definition of $u^0$ there exists a sequence $k^n\in\cK^0$ such that $Y^{0,{\bf 0},k^n}_0\uparrow u^0(0,{\bf 0})$ as $n$ goes to infinity, and  $Y^{0,{\bf 0},k^n}_0\geq u^0_0-c/4$ for all $n$. Define the optimal stopping time $D^k$ for $Y^{0,{\bf 0},k}$ by $D^k=\inf \{s\geq 0 :Y^{0,{\bf 0},k}_s=h_s\}\wedge T$. We can write 
\beaa
Y^{0,{\bf 0},k}_0=\dbE^{\dbP^{0,{\bf 0},k}}\left[\int_0^{D^k\wedge \ch_\d}F^k_r(Y^{0,{\bf 0},k}_r,Z^{0,{\bf 0},k}_r)dr+h_{D^k}\1_{\{D^k< \ch_\d\}}+u^0_{\ch_\d}\1_{\{D^k\geq \ch_\d\}} \right].
\eeaa
Using the uniform bounds on $(Y^{0,{\bf 0},k},Z^{0,{\bf 0},k})$, we have that $$\dbE^{\dbP^{0,{\bf 0},k}}\left[\int_0^{D^k\wedge \ch_\d}|F^k_r(Y^{0,{\bf 0},k}_r,Z^{0,{\bf 0},k}_r)|dr\right]\leq C\sqrt{\d}.$$
uniformly in $k$. We choose $\d$ small such that the previous term is dominated by $\frac{c}{4}$. Recall also that $h_{D^{k_n}}\leq u^0_0-\frac{c}{4}$ on $\{D^{k_n}<\ch_\d\}$ and $u^0_{\ch_\d}\geq u^0_0-\frac{c}{4}$. 

Then, for all $n$, $\dbP^{0,{\bf 0},k_n}(D^{k_n}<\ch_\d)=0$. Therefore $K^{0,{\bf 0},k^n}_{\ch_\d}=0,~\dbP^{0,{\bf 0},k^n}-$a.s. for all $n$. Injecting this into the previous inequalities, the following holds under $\dbP^{0,{\bf 0},k^n}$: 
\beaa
&&(\phi-u^0)_{\ch_\d}+\frac{c \ch_\d}{2}\leq (u^0-Y^{0,{\bf 0},k_n})_{0}+\int_0^{\ch_d}(\d Z^{k_n})^*_s ( dB^t_s+\si_s({k^n_s}) \a_s ds).
\eeaa
There exists a probability $\dbP^n\in\cP_{L_0}$ equivalent to $\dbP^{0,{\bf 0},k^n}$ such that the previous integral is a $\dbP^n$ martingale. Taking the expectation $\dbE^n$:
$$\dbE^n [(\phi-u^0)_{\ch_\d}]+\dbE^n[\frac{c \ch_\d}{2}] \leq (u^0-Y^{0,{\bf 0},k^n})_{0} $$
$\phi\in\underline\cA^{L_0}(0,{\bf 0})$ implies that $0\leq \underline\cE^{L_0}\big[(\phi- u^{0})_{\ch_\d}\big]\leq \dbE^{n} [(\phi-u^0)_{\ch_\d}] $. 
Therefore : $$0<\frac{c}{2}\underline \cE^{L_0}[\ch_\d]\leq \dbE^{n}[\frac{c \ch_\d}{2}] \leq (u^0-Y^{0,{\bf 0},k^n})_{0}. $$
Taking the limit as $n$ goes to infinity we arrive to the contradiction $0<\frac{c}{2}\underline \cE^{L_0}[\ch_\d]\leq0$. 

\subsubsection{Supersolution property of the value functional}

Without loss of generality, we can assume that $$F\mbox{ and }G \mbox{ are decreasing in }y.$$

 We will again reason by contradiction. Assume that $u^0$ is not a viscosity supersolution. A fortiori, it is not a viscosity $L_0$-supersolutions. So there exist $(t,\o)\in[0,T)\times\O$, and $\phi\in\overline \cA^{L_0}u^0 (t,\o)$ with the associated  $\ch\in\cH^t$ such that $-c:=\min (\cL^{t,\o} \phi (t,0),\phi (t,0)-h(t,\o))<0$. Notice that $0=\phi(t,0)-u^0 (t,\o)$ so $\phi(t,0)\geq h(t,\o)$. Therefore $-c=\cL^{t,\o} \phi (t,0)$. 

Without loss of generality we assume $(t,\o)=(0,{\bf 0})$. Similarly to the previous case there exist $\d>0$, such that for $s\in[0,\ch_\d]$ it holds that : 
$$\pa_t \phi_s+ G_s(u^0,\pa_\o\phi)\geq c/2.$$
By the definition of $G$ (in \reff{defn_generator}) and the right continuity of the processes involved, there exists a constant process $k^0\in \cK^0$ such that by taking $\d>0$ small enough, for all $s\in [0,\ch_\d]$ the following inequality holds : 
$$\pa_t \phi_s+\frac12 \pa_{\o\o} \phi_s:\si_s(k^0_s)^2+F^{k^0}_s(u^0_s,\si_s(k^0_s)\pa_\o\phi_s)\geq c/3.$$
We use \reff{second_dpp} with $\t=\ch_\d$ and denote $$(\cY,\cZ,\cK):=(\cY^{0,{\bf 0},k^0}(\ch_\d,u^0_{\ch_\d}),\cZ^{0,{\bf 0},k^0}(\ch_\d,u^0_{\ch_\d}),\cK^{0,{\bf 0},k^0}(\ch_\d,u^0_{\ch_\d}))$$ and with the obvious modifications of the notations of the subsolution case, under $\dbP^{0,{\bf 0},k^0}$ we have:
\beaa
d (\phi- \cY)_s&=&  \Big[(\pa_t \phi_s+\frac12 \pa_{\o\o}\phi_s :\si_s(k^0_s)^2+F_s^{k^0}(u^0_s,\si_s(k^0_s)\pa_\o\phi_s)\Big]ds\notag\\
&&+\Big[F_s^{k^0}(\cY_s,\cZ_s)-F_s^{k^0}(u^0_s,\si_s(k^0_s)\pa_\o\phi_s)\Big]ds+ (\d Z^{k^0}_s)^* dB_s +d\cK_s^{{0,{\bf 0},k^0}}\nonumber\notag\\
&&\geq \frac{c}{6}ds +\Big[F_s^{k^0}(u^0_s,\cZ_s)-F_s^{k^0}(u^0_s,\si_s(k^0_s)\pa_\o\phi_s)\Big]ds+ (\d Z^{k^0}_s)^* dB_s\\
&&\geq \frac{c}{6}ds +(\d Z^{k^0}_s)^*( dB_s+\si_s(k^0_s)\a_s ds)
\eeaa
for some $|\a_s|\leq L_0$.
Therefore under $\dbP^{0,{\bf 0},k^0}$:  
\beaa
&&(\phi-u^0)_{\ch^\d}- (u^0-\cY)_{0}\\
&&=(\phi-u^0)_{\ch^\d}+(u^0-\cY)_{\ch_\d}- (u^0-\cY)_{0}\\
&&=(\phi-\cY)_{\ch_\d}- (u^0-\cY)_{0}\\
&&=(\phi-\cY)_{\ch_\d}- (\phi-\cY)_{0}\geq\frac{c \ch_\d}{6}+\int_0^{\ch_\d} (\cZ^{k^0}_s)^* (dB_s-\si_s(k^0_s)\a_s ds) 
\eeaa
Recall that the DPP \reff{second_dpp} gives $(u^0-\cY)_{0} \geq 0 $ therefore : 
\beaa
&&(\phi-u^0)_{\ch^\d}\geq \frac{c \ch_\d}{6}+\int_0^{\ch_\d} (\cZ^{k^0}_s)^* (dB_s-\si_s(k^0_s)\a_s ds) 
\eeaa
Similarly to the subsolution case, there exist $\dbP\in\cP_{L_0}$, equivalent to $\dbP^{0,{\bf 0},k^0}$ such that the last integral is a $\dbP$ martingale and by assumption $\phi\in\overline\cA^{L_0}u^0(0,{\bf 0}) $. Therefore  $$0\geq \overline \cE^{L_0}[(\phi-u^0)_\ch]\geq \dbE^{\dbP} [(\phi-u^0)_\ch] \geq \frac{c }{6}\dbE^\dbP[\ch]\geq \frac{c }{6}\underline \cE^{L_0}[\ch]>0$$
which is impossible.

\qed
\section{Partial comparison}
Following the same method as \cite{ETZ2}, we will first prove a weaker version of the comparison principle, when, one of the functionals is "smoother". Then we will extend it to general sub/supersolution. In our case the 2RBSDE provides us with a representation formula, therefore our set of "smoother" functionals is simpler than the one in \cite{ETZ2}. Another difference comes from our definition of subsolutions that is only required when the functional does not touch the barrier. Except these points the proofs are the same as the ones in \cite{ETZ2}. We define the following classes of "smoother" processes.
\begin{defn}
\label{defn_Cbar12}
Let $u \in\dbL^0(\L)$, we say that $u$ is in $C^{1,2,-} (\L)$(respectively, $C^{1,2,+} (\L)$) if:\\
(i) There exists a sequence of hitting times $\{\ch_i\}_{i\in\dbN}$, such that $0=\ch_0\leq \ch_1\leq\ldots \leq T$ and the set $\overline\cE^{L}[\1_{\{\ch_n<T\}}]\to 0$ as $n\to \infty$.\\
(ii)For all $(t,\o)\in\L$, $t<T$, and $i$ such that $\ch_i(\o)<t< \ch_{i+1} (\o)$, $u^{t,\o} \in C^{1,2} (\L^t(\ch^{t,\o}_{i+1} ))$,
where $\L^t(\ch^{t,\o}_{i+1} ):=\{(s,\tilde \o)\in\L^t: \ch^{t,\o}_{i+1}(\tilde \o)>s \}$.\\
(iii) For all $i$ given and $r>0$ which is small enough there exists a hitting time $\ch_{i-1}\leq \ch_{i}^{r}\leq \ch_i$ such that for all $\omega$ one has $u^{\ch_{i-1}(\o),\o}$ is uniformly continuous on $[\ch_{i-1}(\o),(\ch_{i}^{r})^{\ch_{i-1}(\o),\o}]$ and $$\overline\cE^{L}_{\ch_i(\o)}\left[|( u_{\ch_{i+1-}}-u_{\ch_{i+1}^{r}})^{\ch_i(\o),\o}|\right]\to 0$$ as $r\to 0$. \\
(iv) $u$ is right continuous and has only nonnegative (resp. nonpositive) jumps at $\ch_i$. 
\end{defn}
\begin{rem}
Notice that we do not require regularity at $\ch_{i}$. However we require right continuity. 
\end{rem}
\begin{thm}\label{Partial Comparison}
Let $u_1\in\underline \cU$ and $u_2\in\overline \cU$ be respectively a viscosity subsolution and a supersolution of \reff{PPDE} such that for all $\o\in \O$, $u^1 (T,\o)\le u^2 (T,\o)$. Assume further that $u_1 \in C^{1,2,-}(\L)$ or $u^2 \in C^{1,2,+}(\L)$, then for all $(t,\o)\in \L$ 
$$u^1(t,\o)\leq u^2 (t,\o).$$
\end{thm}
\proof  To avoid repeating same arguments as in \cite{ETZ2}, we will use the same notations is in the proof of partial comparison in \cite{ETZ2}. We will only point out the differences.

Remark \ref{monotonicity_f} allows us, without loss of generality, to assume that $F$ is non-increasing in $y$. We will prove the statement at $(t,\o)=(0,{\bf 0})$, it is also valid for all intermediate $(t,\o)\in\L$. \\
Define $\hat u:=u^1-u^2$ and denote by $\{\ch_i\}_{i\in\dbN}$ the stopping times given by \reff{defn_Cbar12}. We will first prove that $$\hat u^+_{\ch_i}(\o)\leq \overline\cE^{L}_{\ch_i(\o)}\left[(\hat u^+_{\ch_{i+1-}})^{\ch_i(\o),\o}\right].$$
We only prove the inequality for $i=0$, the proof is valid for all $i$. Assume on the contrary that $$\hat u^+_0 -\overline \cE^L_{0}\left[ \hat u^+_{\ch_1-}\right]>0.$$
By $(iii)$ of definition \ref{defn_Cbar12} one can find $r>0$ small enough such that $$2Tc:=\hat u^+_0 -\overline \cE^L_{0}\left[ \hat u^+_{\ch_1^r}\right]>0.$$
We define $X$ by :
\beaa
&&X:\L\rightarrow \dbR\\
&&X(t,\o):=(\hat u)^+ (t,\o) +ct 
\eeaa 
and  (iii) of definition \ref{defn_Cbar12},  allows us to claim that $X\in \underline \cU$ on $ [0,\ch_{1}^{r}]$. Now define $\hat X:=X\1_{[0,\ch_1^r)}+X_{\ch_1^r-}\1_{[\ch_1^r,T]}$, $Y:=\overline \cS^L [\hat X]$, $\t^*:=\inf\{s\geq 0: Y_t=\hat X_t\}.$\\
Similarly as in \cite{ETZ2}, there exists $\o^*\in\O$ such that $t^*=\t^*(\o^*)<\ch_1^r(\o^*)$ and 
\bea\label{diff_positive}
0<(u^1-u^2)_{t^*}^+ (\o^*)=(u^1-u^2)_{t^*} (\o^*) .
\eea
$X^{t^*,\o^*}\in\underline \cU^{t^*}$ therefore there exists $\d>0$ such that $\ch^{t^*}_\d\leq (\ch_1^r)^{t^*,\o^*}$ and for all $s \in [t^*,\ch^{t^*}_\d]$ it holds that $(u^1-u^2)^{t^*,\o^*}\geq 0$. So we can write  $X^{t^*,\o^*}_t=(u^1-u^2)_t^{t^*,\o^*}+ct$ on $[t^*,\ch^{t^*}_\d].$ 

There are 2 cases to treat. \\
$\bullet$ Assume that $u^2 \in  C^{1,2,+}(\L)$. Then, by definition of $ C^{1,2,+}(\L)$, one has $\phi_t:=(u^2)_t^{t^*,\o^*}-ct \in C^{1,2}(\L^{t^*} (\ch_{1}))$. Then  for all  $ \t\in \cT^{t^*}$ we have :
$$(u^1)_{t^*}^{t^*,\o^*}-\phi_{t^*}=Y_{t^*}(\o^*) \geq \overline \cE^{L}_{t^*}[X^{t^*,\o^*}_{ \t\wedge\ch^{t^*}_\d}],$$ 
which shows that $\phi\in\underline \cA u^1(t^*,\o^*)$. Additionally $u^2 (t^*,\o^*)-h(t^*, \o^*)\geq 0$, so the inequality (\ref{diff_positive}) gives that $(u^1)_{t^*}^{t^*,\o^*} >h(t^*,\o^*)$(this point is the only difference between our proof and the proof in \cite{ETZ2}) and by viscosity subsolution property of $u^1$: 
\beaa
&&0\geq-\pa_t \phi (t^*,\o^*)-G(t^*,\o^*,u^1 (t^*,\o^*),\pa_\o \phi (t^*,\o^*),\pa_{\o\o} \phi (t^*,\o^*))\\
&&=\pa_t u^2 (t^*,\o^*)+ c-G(t^*,\o^*,u^1 (t^*,\o^*),\pa_\o u^2 (t^*,\o^*),\pa_{\o\o} u^2 (t^*,\o^*))\\
&&\geq\pa_t u^2 (t^*,\o^*)+c-G(t^*,\o^*,u^2 (t^*,\o^*),\pa_\o u^2 (t^*,\o^*),\pa_{\o\o} u^2 (t^*,\o^*))\\
&&> 0
\eeaa
which is impossible. \\
$\bullet$ Assume that $u^1 \in  C^{1,2,-}(\L)$. This case is the same as the one in \cite{ETZ2}.

In conclusion, $$\hat u^+_{\ch_i}(\o)\leq \overline\cE^{L}_{\ch_i(\o)}\left[(\hat u^+_{\ch_{i+1-}})^{\ch_i,\o}\right].$$
Then by the Lemma 5.2 of \cite{ETZ2}(which only depends on regularity of $u^1,u^2$ and not on their viscosity solution properties), we have that for all $\dbP\in\cP_L$, it holds that $$\dbE^\dbP\left[\hat u ^+_{\ch_i-}\right]\leq \overline \cE^L_0\left[\hat u^+_{\ch_{i+1-}}\right],$$
By taking the supremum in $\dbP$ and taking into account the positive sign of the possible jumps of $\hat u$ we have that $$\hat u^+_0\leq \overline \cE^L_0\left [\hat u^+_{\ch_n-}\right].$$
By assumption, $\hat u$ is bounded and $\hat u_T^+=0$. Thus using (i) of definition \ref{defn_Cbar12} and passing to the limit in $n$ we obtain $$\hat u^+_0\leq \overline \cE^L_0\left [\hat u^+_{T}\right]=0$$ which completes the proof.  
\qed

\section{Stability}
In this section, we will prove an extension of Theorem 5.1 of \cite{ETZ1} to the PPDE \reff{PPDE}.
\begin{thm} Fix $L>0$ and for $\e>0$, let $(G^\e,h^\e,\xi^\e)$ be a family of data verifying assumptions \reff{assumption_data} with the same constants $M_0,L_0$ and $\rho_0$ and $u^\e$ an $L$-subsolution of \reff{PPDE}. Assume that as $\e$ goes to $0$ the following locally uniform convergences hold :
\bea\label{unif_conv}
&&\mbox{for all }(t,\o,y,z,\g)\in \L\times \dbR\times\dbR^d\times \dbS^d,\mbox{there exists }\d \mbox{ such that :}\\ 
&&(G^\e)^{t,\o}\rightarrow G^{t,\o},\quad (h^\e)^{t,\o}\rightarrow h^{t,\o},\quad(\xi^\e)^{t,\o}\rightarrow \xi^{t,\o},\quad (u^\e)^{t,\o} \rightarrow u^{t,\o},\notag\\
&&\mbox{uniformly on }O^\d_{t,\o,y,z,\g}:=\left \{(s,\o',y',z',\g')\in\L^t\times \dbR\times \dbR^d\times \dbS^d\right.:\notag\\
&&\left .d_\infty^t \Big((s,\o'),(t,{\bf 0})\Big)+|y-y'|+|z-z'|+|\g-\g'|\leq \d\right\},\notag
\eea
Then $u$ is a viscosity $L$-subsolution of the PPDE (\ref{PPDE}) with data $(G,h,\xi)$.
\end{thm}
\begin{rem}
{\rm We are not able to prove the stability result when $L$ depends on $\e$.\qed}
\end{rem}
\begin{rem}
{\rm Except the condition $u(t,\o)\geq h(t,\o)$, our definition of viscosity supersolution is the same as the one given in \cite{ETZ2}. Therefore their stability result for viscosity supersolutions can directly be applied for the PPDE \reff{PPDE}.\qed}
\end{rem}
\proof We will use the same notations as in \cite{ETZ1} and only point out the differences. We will prove the viscosity subsolution property at $(0.{\bf 0})$. We assume that $u(0,{\bf 0})>h(0,{\bf 0})$, $\phi\in\overline \cA^L u(0,{\bf 0})$, and $\ch\in\cH$. The main difference with the proof of Theorem 5.1 in \cite{ETZ1} is that we need to take $\e_\d>0$ small enough to have $u^\e_s>h^\e_s$ for $s\in[0,\ch_\d]$ for all $0<\e<\e_\d$.  Then we have $\phi^\e_\d\in \overline \cA^L u^\e(t^*,\o^*)$ and with  our choice of $\e_\d$ the process $u^\e$ does not touch the barrier $h^\e$. Therefore we can use the viscosity subsolution property of $u^\e$ for the PPDE with data $(G^\e,h^\e,\xi^\e)$ to obtain the equation (5.3) in \cite{ETZ1} and conclude. \qed

\section{Comparison}
Our objective in this section is to extend the partial comparison result. We will carry out the proof in a similar way as in \cite{ETZ2}, and for $0\leq t_1<t_2\leq T$, $\zeta\in\L^0(\dbF_{t_2})$ and $\o\in \O$,  define the following sets :
\beaa
&&\overline \cD(t,\o):=\{\phi\in C^{1,2,+}(\L^t);\min\{\cL^{t,\o}\phi_s,\phi_s-h^{t,\o}_s\}\geq 0,s\in[t,T],\phi_{T}\geq \xi^{t,\o} \},\\
&&\underline \cD(t,\o):=\{\psi\in C^{1,2,-}(\L^{t}); \min\{\cL^{t,\o}\psi_s,\psi_s-h^{t,\o}_s\}\leq 0,s\in[t,T],\psi_{T}\leq \xi^{t,\o} \}.
\eeaa
and processes: 
\beaa
&&\overline u (t,\o):=\inf \{\phi(t,{\bf 0 }): \phi\in\overline \cD(t,\o)\},\\
&&\underline u (t,\o):=\sup\{\psi (t,{\bf 0}):\psi \in \underline \cD(t,\o)\}
\eeaa
\begin{lem}\label{equality}
Under Assumptions \reff{assumption_data}, \reff{assumption_data_add}, the equality $\overline u=\underline u$ holds.
\end{lem}
\proof The proof of this lemma is very technical and requires the introduction of various notations. The construction of the smooth approximating subsolutions and supersolutions are the subject of Appendix A. In the Appendix B, we prove the required regularity of these approximating sequences. 
\qed
\begin{thm}\label{comparison}
Assume \reff{assumption_data}, and \reff{assumption_data_add}, and let $u^1\in \underline \cU$ (respectively, $u_2\in \overline \cU$ ) a viscosity subsolution (respectively, supersolution) of (\ref{PPDE}), such that $u^1 (T,\o)\leq \xi(\o)\leq u^2 (T,\o)$ for all $\o\in\O$, then $u^1 (t,\o)\leq u^2 (t,\o)$ for all $ (t,\o)\in \L$. 
\end{thm}
\proof For all $(t,\o)\in \L$, and $\psi,\phi$ belonging respectively to $\underline \cD (t,\o)$ and $\overline \cD(t,\o)$, by partial comparison result, $u^1 (t,\o)\leq \phi (t,{\bf 0})$ and $\psi (t,{\bf 0})\leq u^2 (t,\o)$. We take the supremum in $\psi$ and the infimum in $\phi$ to have $u^1 (t,\o)\leq \overline u (t,\o)$ and $\underline u (t,\o)\leq u^2 (t,\o)$. The lemma \ref{equality} gives the the equality $\overline u (t,\o)=\underline u(t,\o)$, therefore : $$u^1 (t,\o)\leq u^2 (t,\o)$$
\qed
\appendix
\section{Appendix A}
In the following 2 subsections we will construct 2 families of processes $\{\Psi^{m,\a}\}_{\a>0,m\in\dbN}\in\underline \cD(0,{\bf 0})$ and ${\{\Phi^\a\}}_{\a>0}\in\overline \cD(0,{\bf 0})$ that will allow us to show the Lemma \ref{equality}. We will use the following strategy to prove the equality $\underline u=\overline u=u^0$.
We will freeze the data of the problem $(F,h,\xi)$, in regions of $\L$ related to the stopping times $\ch^t_\d$. Then, we will show that the functionals defined as the solutions of the problem with frozen data are stepwise Markovian. This will bring us to a PDE problem. We will show that the solutions of the given PDEs have interior regularity which will show that the constructed process are in $C^{1,2,+}(\L)$ or $C^{1,2,-}(\L)$.

We recall that, for the comparison result, $\si$ does not depend on $(t,\o)$, and the assumptions on the data allows us to claim that 
\bea\label{def_c_0}
c_0:=\inf_{k\in\cK}\inf _{|\xi|=1}\xi^*\si(k)\xi>0
\eea
 and $F$ is uniformly continuous in $(t,\o)$ with modulus $\rho_0$. Additionally, recall that Remark \reff{change_var} allows us to assume without loss of generality that $F,G$ and $h$ verify \reff{monotonicity_f}.

We will need the following definitions to carry out this construction. For $\a>0$(that will go to $0$) and $t\in[0,T)$, we define:
\beaa
&&\cO^\a:=\{x\in\dbR^d : |x|<\a\},\quad \overline \cO^\a :=\{x\in\dbR^d: |x|\leq \a \},\quad \pa \cO^\a:=\{x\in\dbR^d: |x|=\a\}\\
&&\cO^\a_t:=(t, (t+\a)\wedge T)\times \cO^\a, \quad {\overline \cO^\a_t ;= (t, (t+\a)\wedge T]\times \overline \cO^\a}\\
&&\pa \cO_t^\a ;= ((t, (t+\a)\wedge T]\times \pa \cO^\a)\cup (\{(t+\a)\wedge T\}\times \cO^\a).
\eeaa
For $\{t_i\}_{i\geq 0}$ a nondecreasing sequence in $[0,T]$ with $t_0=0$,  and $\{x_i\}_{i\geq 0}$ a sequence in $\overline\cO^\a$ with $x_0=0$ and $n\geq 0$, we denote $\pi_n :=\{(t_i,x_i)\}_{0\leq i \leq n}$.  In the sequel $\pi_n$ will always verify the previous properties. The sequence $\{t_i\}$ will represent the successive hitting times of a given level by the canonical process, and $\{x_i\}$ the direction of variation of the canonical process between the hitting times. For such $\pi_n$, and $(t,x)\in\cO_{t_n}^\a$, we define  :
\bea\label{hitting_times}
&&{\ch^{t,x,\a}_{-1}:= t_n},\\
&&\ch^{t,x,\a}_0:= \inf \{s\geq t : |x+B_s^t|=\a\}\wedge (t_n+\a)\wedge T,\mbox { and for } i\geq 0,\\
&& \ch^{t,x,\a}_{i+1}:=\inf\{s\geq \ch^{t,x,\a}_i :|B_s^t-B_{\ch_i^{t,x,\a}}^t|=\a\}\wedge (t_n+\a)\wedge T.
\eea

The hitting times $\ch_i^{t,x,\a}$ also depend on $t_n$. For ease of notation we choose not to mention this dependence in the notation. Notice that we can associate to $\pi_n$ a path $\hat\pi_n\in \O$, which is the linear interpolation of $(t_i,\sum_{j=0}^i x_j )_{0\leq i\leq n}$ and $(T,\sum_{j=0}^n x_j )$ and we can associated to $\pi_n$, $(t,x)\in \overline\cO^\a_{t_n}$ and a path $\o\in \O^t$ a  path , $ \hat \o^{\pi_n,t,x,\a}\in \O$, the linear interpolation of $(t_i,\sum_{j=0}^i x_j )_{0\leq i\leq n}$ and of $ (\ch^{t,x,\a}_i (\o), \displaystyle \sum_{j=0}^n x_j+x+\o_{\ch_i^{t,x,\a} (\o)})_{i\geq 0}$.   For $(t,x)\in\overline \cO^\a_{t_n}$, the notation $\pi_n^{(t,x)}$ means that we add $(t,x)$ to the sequence $\pi_n$ as $(n+1)^{th}$ element, namely $\pi_n^{(t,x)}=\{\pi_n, (t,x)\}$. 
For $\pi_n$, and $(t,x)\in\overline \cO^\a_{t_n}$, we define the following generator, final condition and barrier for the approximated equations :  
\beaa
&&\hat F^{\pi_n,t,x,\a}:\L^t\times \dbR\times\dbR^d\times \cK\rightarrow \dbR,\\
&&\hat h^{\pi_n,t,x,\a}:\L^t\times \dbR\rightarrow \dbR,\\
&&\hat \xi^{\pi_n,t,x,\a}:\O^t\rightarrow \dbR.
\eeaa
If  $(s,\o)\in\L^t$, with $\ch^{t,x,\a}_i (\o) \leq s <\ch^{t,x,\a}_{i+1} (\o)$ :
\bea\label{def_hat_f}
&&\hat F^{\pi_n,t,x,\a}(s,\o,y,z,k)=F(\ch^{t,x,\a}_i (\o), \hat \o^{\pi_n,t,x,\a} ,y,z,k), \\
&&\hat h^{\pi_n,t,x,\a}(s,\o)=h(\ch^{t,x,\a}_i (\o), \hat \o^{\pi_n,t,x,\a}), \\
&&\hat \xi^{\pi_n,t,x,\a}(\o)=\xi( \hat \o^{\pi_n,t,x,\a}),
\eea
We remark that $\hat \o^{\pi_n,t,x,\a}$ is not adapted to the filtration $(\dbF^t_s)_{t\leq s\leq T}$. Indeed to know the value of $\hat \o^{\pi_n,t,x,\a}$ after the date $\ch_i^{t,x,\a} (\o)$ we need to know the value of $\o$ at the date $\ch_{i+1}^{t,x,\a} (\o)$. However the data in \reff{def_hat_f} are adapted.

We list the 3 important featured of this approximation:\\
$\bullet$ The approximated generator and barrier are still adapted to $\dbF$.\\
$\bullet$ They verify the assumptions of \cite{Hamadene}, therefore we can use the results on RBSDEs.\\
$\bullet$ Their difference from the original data is less than $\rho_0 (2\a)$.

The idea which consists in approximating the data and studying the RBSDE with the approximated data can not allow us to construct a sequence of subsolutions in $\underline \cD(t,\o)$. Indeed the barrier for the approximated problem might have negative jumps therefore the construction produce subsolutions which might not be in $\underline \cU$ and we would not be able to use this sequence in partial comparison. However this idea allows to produce supersolutions. Indeed the solutions of the RBSDEs can only have negative jumsp which does not create any problem for supersolutions in partial comparison. 

\subsection{Construction of subsolutions by penalization}
In this subsection, we will construct $\{\Psi^{m,\a}\}_{m>0,\a>0}\in\underline \cD(0,{\bf 0})$. The construction will be done by penalization.

Fix $\a>0$, $n,m\in\dbN-\{0\}$, $\pi_n$ as previously, $(t,x)\in\overline \cO^\a_{t_n}$, $k\in \cK^t$, and define $\dbP^{t,k}$ as follows:
\beaa
&&dX^{t,k}_s=\si(k_s)dB^t_s, \mbox{ under }\dbP^t_0\\
&&X^{t,k}_t=0,\\
&&\mbox{and }\dbP^{t,k}:=\dbP^t_0 \circ (X^{t,k})^{-1}.
\eeaa 

Consider $(\cY_s^{\pi_n,t,x,\a,k,m},\cZ_s^{\pi_n,t,x,\a,k,m})_{s\in[t,T]}$(denoted $(\cY_s,\cZ_s)$ for simplicity), the solution of the following BSDE under $\dbP^{t,k}$:
\bea\label{BSDE-}
&&\cY_s =\xi^{\pi_n,t,x,\a} ( B^t)-\int_s^T \cZ_r^* \sigma^{-1}(\tilde k_r)dB^t_r\\
&&+\int_s^T {\hat F^{\pi_n,t,x,\a}_r(\cY_r,\cZ_r,\tilde k_r) +m (\cY_r-\hat h^{\pi_n,t,x,\a}_r)^- dr}\notag
\eea
(see section \ref{sect-intro-value} for the subtlety on $\tilde k$.)and define $$\theta^{\a,m}_n (\pi_n; t,x):= \sup_{k\in\cK^t}\cY_t^{\pi_n,t,x,\a,k,m}.$$ 
Notice that, a priori, we don't know anything on the regularity of this function $\theta_n^{m,\a}$. We also notice that if $(t,x)\in\pa\cO^\a_{t_n}$ then $\ch_0^{t,x,\a}=t$. Therefore $\hat \o^{\pi_n,t,x,\a}=\hat \o^{\pi_n^{(t,x)},t,0,\a}$ for all $\o\in \O^t$ and $(t,x)\in\pa\cO^\a_{t_n}$, which implies the equality of the data defining $(\cY^{\pi_n,t,x,\a,k,m},\cZ^{\pi_n,t,x,\a,k,m})$ and $(\cY^{\pi_n^{(t,x)},t,0,\a,k,m},\cZ^{\pi_n^{(t,x)},t,0,\a,k,m})$. Therefore :
\bea\label{boundary_subsolution}
\theta^{\a,m}_{n} (\pi_n;{t,x})=\theta^{\a,m}_{n+1} (\pi_n^{(t,x)};t,0),\mbox{ for all }(t,x)\in \pa \cO^\a_{t_n}.
\eea

Our main difficulty in the rest of the paper is the fact that the stopping times $\{\ch^{t,x,\a}_i\}$ does not depend continuously on $(t,x)$. However the regularization effect of the PDE allows us to prove the following result.

\begin{prop}\label{regularity_theta}
For all $\a>0$, $n,m \in\dbN-\{0\}$, $\pi_n$ as previously, the mapping 
\bea\label{boundary-theta}
(t,x)\in\pa\cO^\a_{t_n} \rightarrow \theta^{m,\a}_{n+1}(\pi_n^{(t,x)};t,0)=\theta^{m,\a}_{n}(\pi_n;t,x)
\eea 
is continuous and bounded. Additionally $\theta_n^{m,\a}(\pi_n:\cdot)$ is a $C^{1,2}(\cO^\a_{t_n})$ solution of the PDE 
\bea\label{pde_for_theta}
&&-\pa_t \theta_n^{\a,m}(\pi_n;.) -G(t_n,{\hat\pi_n},\theta_n^{\a,m}(\pi_n;.),\pa_x \theta_n^{\a,m}(\pi_n;.),\pa_{xx} \theta_n^{\a,m}(\pi_n;.))\\
&&-m(\theta_n^{\a,m}(\pi_n;.)-h(t_n,\hat\pi_n))^-=0,\mbox{ for all }(t,x)\in \cO^\a_{t_n},\notag\\
&&\theta^{\a,m}_{n} (\pi_n;{t,x})=\theta^{\a,m}_{n+1} (\pi_n^{(t,x)};t,0),\mbox{ for all }(t,x)\in \pa \cO^\a_{t_n}\notag.
\eea
\end{prop}
\proof The proof this proposition and the proposition \ref{gamma_solution} is in the Appendix B. 
\qed
\begin{rem}
Notice that $\pa\cO^\a_{t_n}\cap\{(t_n,x):|x|=\alpha\}=\emptyset$. Thus $\pa\cO^\a_{t_n}$ is not compact and we cannot claim the uniform continuity of the boundary condition and also the solution. Without any additional assumptions one might not be able to extend $\theta_n^{\a,m}(\pi_n;.)$ to $\{(t_n,x):|x|=\alpha\}$ due to lack of uniform continuity near this set. 
\end{rem}

We can now define the process $\{\Psi^{m,\a}\}$. $\ch_i^{0,0,\a}$ will be denoted by $\ch_i$. For $(s,\o)\in\L$ with ${\ch_n (\o)}\leq s < {\ch_{n+1} (\o)}$, $\pi_n (\o)$ will stand for $\{(\ch_i (\o),\o_{\ch_i (\o)}-\o_{\ch_{i-1} (\o)})\}_{0\leq i\leq n}$, and we define:
$$\Psi^{m,\a}(s,\o)=\theta_n^{\a,m}(\pi_n (\o);s,\o_s-\o_{\ch_n (\o)})-\rho_0(\a).$$ 
We now show that, with the associated sequence of stopping times $\{\ch_i\}$, $\Psi^{m,\e}\in C^{1,2,-}(\L)$. Notice that due to the uniform bound on $\si$ and proposition \ref{regularity_theta} the condition (i) and (ii) of definition \ref{defn_Cbar12} holds. (iv) holds by the continuity in time of the solution of second order BSDEs. To prove (iii) we fix $r\in(0,\a/2]$ and define 
\beaa
\ch^r_{i+1}:=\inf\{s\geq \ch_i :|B_s^t-B_{\ch_i}^t|=\a-r\}\wedge (\ch_i+\a-r)\wedge T
\eeaa
Without loss of generality we prove the convergence of the nonlinear expectation for $i=0$. For all $\d>0$ one has 
\beaa 
\overline \cE^L\left[|\Psi^{m,\a}_{\ch_1}-\Psi^{m,\a}_{\ch_1^r}|\right]\leq C \overline \cE^L[\1_{\{\ch_1^{\a/2}\leq \d\}}]+\overline \cE^L\left[|\Psi^{m,\a}_{\ch_1}-\Psi^{m,\a}_{\ch_1^r}|\1_{\{\ch_1^{r}\leq \d\}}\right].
\eeaa
By proposition \ref{regularity_theta} on $\{\ch_1^{r}\leq \d\}$ the function $\Psi^{m,\a}$ is uniformly continuous. Therefore the last term goes to $0$ as $r\downarrow 0$. Which shows that  $\overline \cE^L\left[|\Psi^{m,\a}_{\ch_1}-\Psi^{m,\a}_{\ch_1^r}|\right]\to 0 $ as $r\to 0$. 
 
The definition of $\{\Psi^{m,\a}\}$ gives that if $(t,\o)\in\L$ with $\o_{\ch_n (\o)}<s < \o_{\ch_{n+1} (\o)}$, then denoting  
\beaa
&&P:=\pa_{t}\theta_n^{\a,m}(\pi_n (\o);t,\o_t-\o_{\ch_n (\o)})\\
&&Q:=\theta_n^{\a,m}(\pi_n (\o);t,\o_t-\o_{\ch_n (\o)})\\
&&R:=\pa_x\theta_n^{\a,m}(\pi_n (\o);t,\o_t-\o_{\ch_n (\o)})\\
&&S:=\pa_{xx}\theta_n^{\a,m}(\pi_n (\o);t,\o_t-\o_{\ch_n (\o)})
\eeaa
one has
\beaa
&&\min\{-\pa_t \Psi^{m,\a}(t,\o) -G(t,\o,\Psi^{m,\a}(t,\o),\pa_\o \Psi^{m,\a}(t,\o),\pa_{\o\o} \Psi^{m,\a}(t,\o)),\\
&&\Psi^{m,\a} (t,\o)-h(t,\o)\}\\
&&=\min\{-P -G(t,\o,Q-\rho_0(\a),R,S),Q-\rho_0 (\a)-h(t,\o)\}\\
&&\leq \min\{-P-G(t,\o,Q, R,S)-\rho_0 (\a),Q-h(\ch_n (\o),\hat \pi_n(\o))\}\\
&&\leq \min\{-P-G(\ch_n (\o),\pi_n(\o),Q,R,S),Q-h(\ch_n (\o),\hat \pi_n(\o))\}\leq 0.
\eeaa
 which shows that $\Psi^{m,\a}\in \underline \cD(0,{\bf 0})$.

\subsection{Construction of supersolutions by approximation}
Fix $\a>0$, $\pi_n$, for $(t,x)\in\overline\cO^\a_{t_n}$ and $k\in\cK^t$, our approximated data defined at (\ref{hitting_times}) verifies the assumption of \cite{Hamadene}, therefore we have the existence of $(\hat Y^{\pi_n,t,x,\a,k},\hat Z^{\pi_n,t,x,\a,k}, \hat K^{\pi_n,t,x,\a,k})_{s\in[t,T]}$ (we drop the superscript ${\pi_n,t,x,\a,k}$ for simplicity of notation) solution of the following RBSDE under $\dbP^{t,k}$. 
\bea\label{rbsde_approximation}
&&\hat Y_s =\hat\xi^{\pi_n,t,x,\a}(B^t) +\int_s^T \hat F^{\pi_n,t,x,\a}_r(\hat Y_r,\hat Z_r,\tilde k_r)dr-\int_s^T \hat Z_r^* \sigma^{-1}(\tilde k_r)dB^t_r +\hat K_T-\hat K_s, \notag\\
&&\hat Y_s\ge \hat h^{\pi_n,t,x,\a}_s,\quad [\hat Y_s-\hat h^{\pi_n,t,x,\a}_s] d\hat K^c_s=0, \\
&& \Delta_s \hat Y:= \hat Y_s-\hat Y_{s-}=-(\hat h^{\pi_n,t,x,\a}_{s-}-\hat Y_s)^+,\quad \hat K  \mbox{ non decreasing, }\hat K_t=0.\notag
\eea
And similarly we define the following mapping : 
$$\Gamma_n^\a (\pi_n;t,x):= \sup_{k\in\cK^t}\hat Y_t^{\pi_n,t,x,\a,k}.$$
Notice that if $\Delta_s \hat Y>0$ than $\Delta_s \hat h^{\pi_n,t,x,\a}>0$. Therefore the jumps of $\hat Y$, which are the jumps of the discontinuous part $\hat K^d$ of $\hat K$, can only happen when there is a jump of $\hat h^{\pi_n,t,x,\a}$, those possible jump dates are $\{\ch^{t,x,\a}_i\}$. \\
In the literature, there are some estimates of $d\hat K^c$, the continuous part of $\hat K$ when the barrier is a continuous semimartingale. In that case it can be shown that $0\leq d\hat K^c_s\leq (\hat F^{\pi_n,t,x,\a}_s(\hat Y_s,\hat Z_s,\tilde k_s) ds + d A_s)^-$ where $A$ is the drift part of the barrier,(notice that in our case $dA_s=0$, excepts at $\{\ch^{t,x,\a}_i\}$). We will extend this result to our case. 

\subsubsection {Study of $K^c$}\label{study_K}
At this subsection, we will study the RBSDE defined at (\ref{rbsde_approximation}) under $\dbP^{t,k}$. We again drop the superscript ${\pi_n,t,x,\a,k}$ for notational simplicity. We denote by $(\hat F, \hat \xi,\hat h)$ the data and by $(\hat Y, \hat Z, \hat K)$ the solution. Recall that by the Remark \ref{rem_change}, we can assume that  $\hat F$ and $\hat h$ verifies (\ref{monotonicity_f}).
We then have the following proposition :
\begin{prop}\label{prop-kc}
Under assumptions \reff{assumption_data},
$$\hat K^c\equiv 0,\q \dbP^{t,k} \mbox{-a.s}.$$
\end{prop}
\proof
%The following process are caglad therefore predictable :
%\bea
%&&\hat h^n(s,\o)=\sum_{i=1}^\infty {\1_{(\ch^n_i(\o),\ch^n_{i+1}(\o)]}(s)h(\ch^n_i(\o),\o)} \mbox { if } t\leq s\leq T-\frac{\e}{n}\\
%&&\hat h^n(s,\o)=\xi(\o{.\wedge T-\frac{\e}{n}}) \mbox { if }  T-\frac{\e}{n} < s\leq T
%\eea
%as $n$ goes to infinity 
%\bea
%&&\1_{(\ch^n_i(\o),\ch^n_{i+1}(\o)]}(s) \rightarrow \1_{[\ch_i(\o),\ch_{i+1}(\o))}(s)\\
%&&h(\ch^n_i(\o),\o)\rightarrow h(\ch_i(\o),\o) 
%\eea
 %which shows that $\hat h^n\rightarrow \hat h$. So $\hat h $ is predictable. It also has finite variation. 
We differentiate $(\hat Y-\hat h)$ in 2 different ways under $\dbP^{t,k}$: 
\bea\label{compute-local}
d(\hat Y_s-\hat h_s)=-\hat F_s ds -d\hat K^c_s - d\hat K^d_s -d\hat h_s +\hat Z_r^* \sigma^{-1}(\tilde k_r)dB^t_r,
\eea
where $\hat F_s := \hat F^{\pi_n,t,x,\a}_s(\hat Y_s,\hat Z_s,\tilde k_s)$, and $\hat h_s=\hat h^{\pi_n,t,x,\a}(s,B^t)$ . The processes $\hat K^c$ and $\hat K^d$ have integrable variation. Between two successive $\ch^{t,x,\a}_i$ (denoted $\ch_i$ for simplicity), $\hat h$ is constant, therefore the variation of $\hat h (\o)$ is bounded by $2(N(\o)+1)\rho_0(\e)<\infty$ for all $\o\in\o^t$, where $N(\o)=\inf\{i\in \dbN: \ch_i(\o)=T\}<\infty$. Additionally $\hat h$ is cadlag and constant between the terms of $\{\ch_i\}$, so $(\hat Y-\hat h)$ is a semimartingale, denoting $L^0$ its local time at $0$, by the It\^o-Meyer formula :
\beaa
&&d(\hat Y_s - \hat h_s)=d(\hat Y_s - \hat h_s)^+=\\
&& =\1_{\hat Y_{s-}>\hat h_{s-}}\big(-\hat F_s ds -d\hat K^d_s -d\hat h_s  +\hat Z_s^* \sigma^{-1}(\tilde k_s)dB^t_s  \big)\\
&&+\1_{\hat Y_{s-} =\hat h_{s-} } \Delta_s (\hat Y-\hat h)+\frac12 L^0_s.
\eeaa
In the previous equality, we used $0=({\hat Y_{s-}-\hat h_{s-}})d\hat K^c_s=\1_{\{\hat Y_{s}>\hat h_{s}\}}d\hat K^c_s$ to eliminate the term $\1_{\hat Y_{s-}>\hat h_{s-}}d\hat K^c_s$.  Identifying with \reff{compute-local} :
\beaa
&&\1_{\hat Y_{s-} =\hat h_{s-} } \big( - \hat F_s dt - d\hat K^d_s -d\hat h_s +\hat Z_s^* \sigma^{-1}(\tilde k_s)dB^t_s  \big) -d\hat K^c_s=\1_{\hat Y_{s-}=\hat h_{s-} } \Delta_s (\hat Y_s-\hat h_s )+\frac12 L_s^0.
\eeaa
We define the set $\cJ:=\{(s,\o): s\neq \ch_i (\o)\mbox { for all } i \}$ and notice that on $\cJ$, $d\hat h_s= d\hat K_s^d =\Delta_s (\hat Y-\hat h)=0$, so $\hat Y_s=\hat Y_{s-}$ and $\hat h_s=\hat h_{s-}$. By rewriting the previous equality on $\cJ$ :
\bea\label{meyer_ito}
&&\1_{\hat Y_{s} =\hat h_{s} } \big( \hat Z_s^* \sigma^{-1}(\tilde k_s)dB^t_s  \big) =d\hat K^c_s+\1_{\hat Y_{s} =\hat h_{s} } \big(   \hat F_s ds  \big)+\frac12 L_s^0.
\eea
The right term is predictable finite variation and the left term defines a martingale. Therefore on the set $\{\hat Y_s=\hat h_s\} \cap \cJ$, we have : $$\hat Z_s=0,\mbox{ and}$$ $$0\leq d\hat K_s^c \leq \hat F_s ds =\hat F_s^-(\hat h_s ,\hat Z_s,\tilde k_s)ds=\hat F^-_s (\hat h_s,0,\tilde k_s)ds.$$
Notice that for $\o\in\O^t$ using  the Remark \ref{rem_change}:
\beaa
&&\hat F_s(\hat h_s,0,\tilde k_s)= \hat F^{\pi_n,t,x,\a}_s (\hat h^{\pi_n,t,x,\a}(s,\o),0,\tilde k_s=F(s,\hat\o^{\pi_n,t,x,\a},h(s,\hat \o^{\pi_n,t,x,\a}),0,\tilde k_s)\geq 0.
\eeaa
 So $d\hat K_s^c=0$ on $\{\hat Y_s=\hat h_s\} \cap \cJ$. 
On  $\{\hat Y_s \neq \hat h_s\} \cap \cJ$, the equality (\ref{meyer_ito}) directly gives $d\hat K_s^c=0$. \\
Therefore on $\cJ$, $d\hat K_s^c=0$, $dt\times\dbP^{t,k}$-a.s. which shows that $\hat K^c$ is constant between the $\ch_i$ so it is always $0$. 
\qed\\
Then $\hat K=\hat K^d$ can only jump at the stopping times $\{\ch_i^{t,x,\a}\}$.

If we rewrite \reff{rbsde_approximation} up to $\ch_0^{t,x,\a}$ for $s<\ch_0^{t,x,\a}$ it becomes (without the superscript $(\pi_n,t,x,\a,k)$) under $\dbP^{t,k}$:
\beaa
&&\hat Y_s =\hat Y_{\ch_0^{t,x,\a}}+\hat K_{\ch_0^{t,x,\a}} +\int_s^{\ch_0^{t,x,\a}} \hat F_r(\hat Y_r,\hat Z_r,\tilde k_r)dr-\int_s^{\ch_0^{t,x,\a}} \hat Z_r^* \si^{-1}(\tilde k_r) dB^t_r,  \\
&&\hat Y_s\ge \hat h_s, \\
&& \hat K_{\ch_0^{t,x,\a}}= - \Delta_{\ch_0^{t,x,\a}} \hat Y:= -(\hat Y_{\ch_0^{t,x,\a}}-\hat Y_{{\ch_0^{t,x,\a}-}})=(\hat h_{{\ch_0^{t,x,\a}}-}-\hat Y_{\ch_0^{t,x,\a}})^+. 
\eeaa
Therefore  for $s<\ch_0^{t,x,\a}$:
\bea\label{actual-bsde}
\hat Y_s =\max \{\hat Y_{\ch_0^{t,x,\a}};   h(t_n,\hat \pi_n)\}+\int_s^{\ch_0^{t,x,\a}} \hat F_r(\hat Y_r,\hat Z_r,\tilde k_r)dr-\int_s^{\ch_0^{t,x,\a}} \hat Z_r^* \si^{-1}(\tilde k_r)  dB_r^t.
\eea
This equation is actually a BSDE up to $\ch_0^{t,x,\a}$. We need the following results to continue our analysis. 
\begin{prop}\label{gamma_solution} Under the assumptions \reff{assumption_data},  for fixed $n$ and $\a>0$,
the mapping $\cM^\a_n$
\bea\label{final_condition_super}
(t,x)\in \pa\cO^{\a}_{t_n}\rightarrow\max \{\Gamma_{n+1}^\a (\pi_n^{(t,x)};t,0); \hat  h(t_n,\hat \pi_n)\}=:\cM^\a_n(\pi_n;t,x),
\eea 
is continuous and bounded and the function $\Gamma_n^\a(\pi_n;.)$ is a $C^{1,2}(\cO^\a_{t_n})$ solution of the following PDE : 
\bea\label{equation_Phi}
&&-\pa_t \G_n^{\a}(\pi_n;.) -G (t_n,\hat\pi_n,\G_n^{\a}(\pi_n;.),\pa_x \G_n^{\a}(\pi_n;.),\pa_{xx} \G_n^{\a}(\pi_n;.))=0,\mbox{ for all }(t,x)\in \cO^\a_{t_n},\notag\\
&&\G^{\a}_{n} (\pi_n;{(t,x)})=\cM^\a_n(\pi_n;t,x),\mbox{ for all }(t,x)\in \pa \cO^\a_{t_n}.\label{boundary_super_sol}
\eea
and satisfies $$\G_n^{\a}(\pi_n;t,x)\geq  h (t_n,\hat\pi_n) \mbox{ for all }(t,x)\in \cO^\a_{t_n}.$$

\end{prop}
\proof
The proof of this result is the subject of the next Appendix. 
\qed

We define $ \Phi^\a \in  C^{1,2,+}(\L)$ as in the subsolution case. Similarly, $\ch_i$ stands for $\ch_i^{0,0,\a}$, for $(s,\o)\in\L$ with ${\ch_n (\o)}\leq s < {\ch_{n+1} (\o)}$, $\pi_n (\o)$ stands for $\{(\ch_i (\o),\o_{\ch_i (\o)}-\o_{\ch_{i-1} (\o)})\}_{0\leq i\leq n}$, and we define :
\beaa
\Phi^\a (s,\o):=\G_n^\a \big(\pi_n (\o); s,\o_s-\o_{\ch_{n} (\o)}\big)+\rho_0 (\a).
\eeaa
As it is proven for $\Psi^{m,\a}$, using the Remark \ref{change_var}, $\Phi^\a\in \overline \cD(0,{\bf 0})$ .

Finally, we can now prove the lemma \ref{equality}.\\
{\bf Proof of lemma \ref{equality}}
We show the two inequalities separately.\\
Notice that by partial comparison, for any $\phi \in \overline \cD(0,{\bf 0})$ and $\psi \in \underline \cD(0,{\bf 0})$, $\phi(0,{\bf 0})\geq \psi(0,{\bf 0})$, which by taking the infimum in $\phi$ and supremum in $\psi$ shows that $ \overline u(0,{\bf 0})\geq \underline u(0,{\bf 0})$. $\Psi^{m,\a}\in\underline \cD(0,{\bf 0})$ and $\Phi^\a(0,{\bf 0})\in\overline \cD(0,{\bf 0})$, so $$\Psi^{m,\a}(0,{\bf 0})\leq \underline u(0,{\bf 0})\leq \overline u(0,{\bf 0}) \leq \Phi^\a(0,{\bf 0}).$$
Fix $\d>0$ and $\a>0$, then $\Gamma^{\a}_0(\pi_0;0,{\bf 0})$ is the value at $0$ of the 2RBSDE with data $(\hat G^{0,0,\a},\hat h^{0,0,\a},\xi^{0,0\a})$ and $\theta^{m,\a}_0(\pi_0;0,{\bf 0})$ is the value at $0$ of the 2BSDE with generator $$\hat G^{0,0,\a}(s,\o,x,y)-m(y-\hat h^{0,0,\a})^-$$ and final condition $\xi^{0,0,\a}$, the convergence of the solutions of the penalized BSDE to the solution of the RBSDE gives that there exists $m_\a\in\dbN$ such that $\theta^{m_\a,\a}_0(\pi_0;0,{\bf 0})\geq \Gamma^\a_0 (\pi_0;0,{\bf 0})-\d$. We rewrite these inequalities in terms of $\Psi^{m_\a}$ and $\Phi^\a$ to have  $\Psi^{m_\a,\a}(0,{\bf 0})+\rho(\a)+\a\geq \Phi^{\a}(0,{\bf 0})-\rho(\a)-\a-\d.$
By the definition of $\overline u$, and $\underline u$, this gives 
$$\underline u(0,{\bf 0})+\rho_0 (\a)+\a\geq \overline u(0,{\bf  0})-\a-\rho_0(\a)-\d.$$
We take take the limit as $\a,\d$ goes to 0 to have $$\underline u(0,{\bf 0})\geq \overline u(0,{\bf 0}).$$

\qed

\section{Appendix B}\label{Appendix B}
In this section we provide a proof for the regularity of $\G^\a_n$. The result also holds for  $\theta^{m,\a}_n$. In this case the generic constant $C$ also depend on $m$.  

We prove the continuity of the boundary condition at Proposition \ref{gamma_solution} in 2 steps.

{\it Step 1: Regularity in space:}
We first prove a lemma on the dependence of the solutions of the approximated problem on $\pi_n$.
\begin{lem}\label{regularity_in_pi_n}
There exist $C>0$ depending only on $d,L_0,M_0$, and $T$ such that for all $\pi_n=\{(t_i,x_i)\}_{0\leq i\leq n}$, and $\pi'_n=\{(t'_i,x'_i)\}_{0\leq i\leq n}$ with $t_n=t'_n$, 
$$|\G^{\a}_{n} (\pi'_n;t_n,0)-\G^{\a}_{n} (\pi_n;t_n,0)|\leq C\rho_0 (||\hat \pi_n-\hat \pi'_n||_{t_n})$$ 
holds. Thus the function $\cM_n^\a(\pi_n;.)$ in \reff{final_condition_super} is uniformly continuous in $x$ with modulus $C\rho_0$. 
\end{lem}
\proof
$\G^{\a}_{n} (\pi'_n;t_n,0)$ and $\G^{\a}_{n} (\pi_n;t_n,0)$ are defined with the solutions at $t_n$ of RBSDEs with respectively data $$(\hat F^{\pi_{n},t_n,0,\a}(s,\o,y,z,k),\hat h^{\pi_{n},t_n,0,\a}(s,\o),\hat \xi^{\pi_{n},t_n,0,\a}(\o))$$ and $$(\hat F^{\pi'_{n},t_n,0,\a}(s,\o,y,z,k),\hat h^{\pi'_{n},t_n,0,\a}(s,\o),\hat \xi^{\pi'_{n},t_n,0,\a}(\o))$$ and the stopping times $\{\ch_i^{t_n,0,\a}\}$ in \reff{hitting_times} are the same for both of the data. Therefore for all $(s,\o)\in\L^t$, denoting $\d= ||\hat\pi_n-\hat\pi'_n||_{t_n}$, we have :
\beaa
&&|\hat F^{\pi_{n},t_n,0,\a}(s,\o,y,z,k)-\hat F^{\pi'_{n},t_n,0,\a}(s,\o,y,z,k)|\leq \rho_0(\d),\\
&&|\hat h^{\pi_{n},t_n,0,\a}(s,\o)-\hat h^{\pi'_{n},t_n,0,\a}(s,\o)|\leq \rho_0(\d),\\
&&|\hat \xi^{\pi_{n},t_n,0,\a}(\o)-\hat \xi^{\pi'_{n},t_n,0,\a}(\o)|\leq \rho_0(\d),
\eeaa
Given the a priori estimates of RBSDEs  and taking the sup in $k\in \cK^{t_n}$, we have that 
\beaa
&&|\G^{\a}_{n} (\pi_{n};t_n,0)-\G^{\a}_{n} (\pi'_{n};t_n,0)|\leq C\rho_0(\d).
\eeaa
\qed

{\it Step 2: Regularity in time}:
\begin{lem}\label{regularity_in_pi_n}
For all $\pi_n=\{(t_i,x_i)\}_{0\leq i\leq n}$ and $(t,x),(t',x)\in \pa\cO^\a_{t_n}$ with $t_n<t\leq t_n+\a$ then, 
\beaa
 \lim_{\substack{t'\to t\\ t'\in(t_n,t_n+\a)}}\G^{\a}_{n} (\pi_n;t',x)=\G^{\a}_{n} (\pi_n;t,x).
\eeaa
\end{lem}
Before proving this lemma, we finish the proof of Propositions \ref{gamma_solution} and \ref{regularity_theta}. 
Given the continuity of the boundary condition, the PDE's in the propositions \ref{regularity_theta} and \ref{gamma_solution} has the following general form 

\bea\label{general-pde}
-\pa_t v(t,x)-g(v(t,x),\pa_{x}v(t,x),\pa_{xx}v(t,x))=0,\mbox{ for }(t,x)\in \cO^\a_0\\
v(t,x)=h(t,x),\mbox{ for }(t,x)\in\pa \cO^\a_0,
\eea
where $g(y,z,\gamma)=\sup_{k}\{\frac{1}{2}\si^2(k):\gamma+f^k(y,z)\}$ with $f^k$ is Lipschitz continuous and the boundary condition $h$ is continuous and bounded.  Thus for all $r>0$ the PDE \reff{general-pde} satisfies the assumptions of theorem 1.1 (and also remark 1.1(i)) of \cite{TW} for $t\in[r,\a]$. Thus using the interior estimates of theorem 1.1 of \cite{TW} the solution of the PDE \reff{general-pde} is $C^{1,2}(\cO^\a_0)$. 

The rest of the paper contains the proof of Lemma \ref{regularity_in_pi_n}. 

Recall that  $(t,x),(t',x)\in \pa\cO^\a_{t_n}$ with $t_n<t \leq t_n+\a$ and $t_n<t' < t_n+\a$ thus $||\widehat{{\pi_n^{t,x}}}-\widehat{{\pi_n^{t',x}}}||_{t\wedge t'}=\a\frac{|t'-t|}{t\wedge t'-t_n}$. The presence of this denominator is the reason why one cannot obtain uniform continuity of the boundary condition on its whole domain. To show the right continuity, we take $t_n<t<t'$, one can also similarly show the left continuity. Notice that $\ch_0^{t,x,\a}=t$ and define $$\ch:=\inf\{s\geq t: |B^t_s|\geq \a/2\}\wedge t'<\ch_1^{t,x,\a}.$$
Given the proposition \ref{prop-kc}, using the DPP for $\G^{\a}_{n} (\pi_n;\cdot)$ and \reff{actual-bsde} one has
$$\G^{\a}_{n} (\pi_n;t,x)=\sup_{k\in\cK^t} \tilde Y^{\pi_n,t,x,\a,k}_t$$
where $(\tilde Y^{\pi_n,t,x,\a,k}_s,\tilde Z^{\pi_n,t,x,\a,k}_s)$ solves the BSDE with random maturity time $\ch$ under $\dbP^{t,k}$
$$Y_s =\G^{\a}_{n+1} (\pi_n^{(t,x)};\ch,B^t_\ch) +\int_s^\ch \hat F^{\pi_n,t,x,\a}_r(Y_r, Z_r,\tilde k_r)dr-\int_s^\ch Z_r^* \sigma^{-1}(\tilde k_r)dB^t_r .
$$
Then under $\dbP^{t,k}$
\beaa
\tilde Y^{\pi_n,t,x,\a,k}_t -\G^{\a}_{n+1} (\pi_n^{(t',x)};t',0)
&=&\G^{\a}_{n+1} (\pi_n^{(t,x)};\ch,B^t_\ch) -\G^{\a}_{n+1} (\pi_n^{(t',x)};t',0)\\
&+&\int_t^\ch \hat F^{\pi_n,t,x,\a}_r(Y_r, Z_r,\tilde k_r)dr-\int_t^\ch Z_r^* \sigma^{-1}(\tilde k_r)dB^t_r 
\eeaa
Taking the expectation under $\dbP^{t,k}$ 
\beaa
|\tilde Y^{\pi_n,t,x,\a,k}_t -\G^{\a}_{n+1} (\pi_n^{(t',x)};t',0)|
\leq\dbE^{\dbP^{t,k}}\left[|\G^{\a}_{n+1} (\pi_n^{(t,x)};\ch,B^t_\ch) -\G^{\a}_{n+1} (\pi_n^{(t',x)};t',0)|\right]+\rho(\sqrt{|t-t'|}) 
\eeaa
We then take $k_n$ such that $\tilde Y^{\pi_n,t,x,\a,k_n}_t\to\G^{\a}_{n+1} (\pi_n;t,x)=\G^{\a}_{n+1} (\pi_n^{(t,x)};t,0)$ to obtain that 
\beaa
|\G^{\a}_{n+1} (\pi_n^{(t,x)};t,0) -\G^{\a}_{n+1} (\pi_n^{(t',x)};t',0)|
&&\leq\sup_{k\in\cK^t}\dbE^{\dbP^{t,k}}\left[|\G^{\a}_{n+1} (\pi_n^{(t,x)};\ch,B^t_\ch) -\G^{\a}_{n+1} (\pi_n^{(t',x)};t',0)|\right]\\&&+\rho(\sqrt{|t-t'|}) .
\eeaa
Notice that by simple estimates $\sup_{k\in\cK^t}{\dbP^{t,k}}(\ch<t')\to 0$ as $t'\to t$ and $\overline\cE^{L}[\1_{\{\ch_N^{t',0,\a}<T\}}]\to 0$ as $N\to \infty$, thus one can find $N$ independent of $t'$ such that 
\bea\label{estimate_gamma_prime}
&&|\G^{\a}_{n+1} (\pi_n^{(t,x)};t,0) -\G^{\a}_{n+1} (\pi_n^{(t',x)};t',0)|\\
&&\leq\sup_{k\in\cK^t}\dbE^{\dbP^{t,k}}\left[\1_{\{\ch=t'\}\cap\{\ch_N^{t'0,\a}=T\}}|\G^{\a}_{n+1} (\pi_n^{(t,x)};t',B^t_{t'}) -\G^{\a}_{n+1} (\pi_n^{(t',x)};t',0)|\right]\notag\\&&+\rho(\sqrt{|t-t'|}) \notag
\eea
The next step is to estimates $|\G^{\a}_{n+1} (\pi_n^{(t,x)};t',B^t_{t'}) -\G^{\a}_{n+1} (\pi_n^{(t',x)};t',0)|$ on ${\{\ch=t'\}\cap\{\ch_N^{t'0,\a}=T\}}$. Notice that on this event $|B^t_{t'}|\leq \a/2$ and for all $y$ such that $|y|\leq \a/2$ and $(s,\o)\in\L^t$ one has that 
\bea\label{estim-bsde}
|\hat F^{\pi_{n}^{(t,x)},t,y,\a}(s,\o,y,z,k)-\hat F^{\pi_{n}^{(t',x)},t',0,\a}(s,\o,y,z,k)|\leq \a\frac{|t'-t|}{t-t_n}+||\hat \o^{t',y,\a}-\hat \o^{t',0,\a}||^{t'}_T.
\eea
Notice that the same inequality holds also for $\hat h$ and $\hat\xi$. We need the next Lemma to control the last term in this inequality by controling the difference $|\ch^{{t'},y,\a}_i-\ch^{{t'},{\bf 0},\a}_i|$.  
\begin{prop}\label{estimates_ch}
For $n>0$ define 
\bea\label{definition_gamma} 
&&\Delta^{{t'},y,\a}_n:=\sup_{0\leq i\leq n}|\ch^{{t'},y,\a}_i-\ch^{{t'},{\bf 0},\a}_i| \mbox{ and}\\
&&\dbC^{t'}_{K,1/3}:=\{\o\in\O^t, \sup_{{t'}\leq s<r\leq T}\frac{|\o_r-\o_s|}{|r-s|^{1/3}}<K  \}
\eea
then for all  $\e>0$, $\d>0$ and $n\in \dbN$ there exist $K_\e <\infty$(depending only on $\e$) and $q>0$ (independent of ${t'}$), such that $\sup_{k\in\cK^{t'}}\dbP^{t',k}(\{\Delta^{t',y,\a}_n>\d\} \cup (\dbC^{t'}_{K_\e,1/3})^c)\leq \e$ if $|y|\leq q$.
\end{prop}
\proof By classical results on stochastic analysis, there exist a constant $p>0$ depending only on $d$ such that $\sup_{k\in\cK^{t'}}\dbE^{\dbP^{t',k}} \left[\sup_{{t'}\leq s<r\leq T}\frac{|B^{t'}_r-B^{t'}_s|^p}{|r-s|^{p/3}}\right]<\infty$. Then $\sup_{k\in\cK^{t'}}{\dbP^{t',k}}((\dbC^{t'}_{K,1/3})^c)\leq \frac{1}{K^p}\sup_{k\in\cK^{t'}}\dbE^{\dbP^{t',k}} \left[\sup_{{t'}\leq s<r\leq T} \frac{|B^{t'}_r-B^{t'}_s|^p}{|r-s|^{p/3}}\right] \rightarrow 0 $ as $K$ goes to infinity. Remark that the upper bound depends only on $d,L_0$ and $T$. Thus we obtain the existence of $K_\e$.   

Similarly to the proof of Lemma 4.6 in \cite{ETZ2} one can find a constant $C$ depending on $c_0:=\inf_{k\in\cK}\inf _{|\xi|=1}\xi^*\si(k)\xi>0$ such that $\sup_{k\in\cK^{t'}}{\dbP^{t',k}}(|\ch^{t',y,\a}_0-\ch^{t',0,\a}_0|\geq \d)\leq C\frac{|y|}{\sqrt{\d}}$. 

We now fix $k\in\cK^{t'}$ and $\d_1,\ldots \d_n>0$ to be determined and $\d_i\leq \d$ for all $i$. 
\beaa
\dbP^{t',k}(\{\Delta^{t',y,\a}_n>\d\} \cup (\dbC^{t'}_{K_\e,1/3})^c)&=&\dbP^{t',k}(\{\Delta^{t',y,\a}_n>\d\} \cap \dbC^{t'}_{K_\e,1/3})+\dbP^{t',k}( (\dbC^{t'}_{K_\e,1/3})^c)\\
&\leq &\dbP^{t',k}(\{\Delta^{t',y,\a}_n>\d\} \cap \dbC^{t'}_{K_\e,1/3})+\frac{\e}{2}\\
&\leq &\sum_{i=0}^n \dbP^{t',k}(\{|\ch^{{t'},y,\a}_i-\ch^{{t'},{\bf 0},\a}_i|>\d\} \cap \dbC^{t'}_{K_\e,1/3})+\frac{\e}{2}\\
&\leq &\sum_{i=0}^n \dbP^{t',k}(\{|\ch^{{t'},y,\a}_i-\ch^{{t'},{\bf 0},\a}_i|>\d_i\} \cap \dbC^{t'}_{K_\e,1/3})+\frac{\e}{2}
\eeaa
We estimate the term inside the sum 
\beaa
&&\dbP^{t',k}(\{|\ch^{{t'},y,\a}_i-\ch^{{t'},{\bf 0},\a}_i|>\d_i\} \cap \dbC^{t'}_{K_\e,1/3})\\
&&\q\leq \dbP^{t',k}(\{|\ch^{{t'},y,\a}_{i-1}-\ch^{{t'},{\bf 0},\a}_{i-1}|>\d_{i-1}\} \cap \dbC^{t'}_{K_\e,1/3})\\
&&\q+\dbP^{t',k}(\{|\ch^{{t'},y,\a}_{i-1}-\ch^{{t'},{\bf 0},\a}_{i-1}|\leq \d_{i-1}\} \cap \{|\ch^{{t'},y,\a}_i-\ch^{{t'},{\bf 0},\a}_i|>\d_i\} \cap \dbC^{t'}_{K_\e,1/3})
\eeaa
Notice that on $ \dbC^{t'}_{K_\e,1/3} $ we have $|B^{\ch^{t',y,\a}_{i-1}\wedge \ch^{t',0,\a}_{i-1}}_{\ch^{t',y,\a}_{i-1}\vee \ch^{t',0,\a}_{i-1}}|\leq K_\e |\ch^{t',y,\a}_{i-1}- \ch^{t',0,\a}_{i-1}|^{1/3}$. Thus we can control the last term with $\frac{C K_\e \d_{i-1}^{1/3}}{\sqrt {\d_i}}$. 
Thus one can find universal constants $C_{n,i}$ such that 
\beaa
\dbP^{t',k}(\{\Delta^{t',y,\a}_n>\d\} \cup (\dbC^{t'}_{K_\e,1/3})^c)&\leq&\sum_{i=1}^n C_{n,i}\frac{ K_\e \d_{i-1}^{1/3}}{\sqrt {\d_i}}+C\frac{|y|}{\sqrt{\d_0}}
\eeaa
We now choose $\d_n=\d$ and inductively determine $\d_i$, $i=n-1,\ldots, 0$ small enough to have the previous sum small enough. After determining $\d_0$ we determine $q>0$ as required.
 \qed
 
We continue with the proof of Lemma \ref{regularity_in_pi_n}.
Using classical estimates on BSDEs and \reff{estim-bsde}(for $\hat F$, $\hat h$ and $\hat \xi$) one has
\bea\label{estimate-finall}
|\G^{\a}_{n+1} (\pi_n^{(t,x)};t',y) -\G^{\a}_{n+1} (\pi_n^{(t',x)};t',0)|&\leq&\frac{C\a |t-t'|}{t-t_n}+C\overline\cE^{L}[\1_{\{\ch_N^{t',0,\a}<T\}}]\\
&+&C\sup_{k\in\cK^{t'}}\dbE^{\dbP^{t',k}}\left[(||\hat \o^{t',y,\a}-\hat \o^{t',0,\a}||^{t'}_T)^2\1_{\{\ch_N^{t',0,\a}=T\}}\right]\notag
\eea
We now fix $\e>0$ and determine $K_\e$ as in the previous proposition. We then choose $\d=\e^3/K_\e^3$ and apply again the previous lemma for $n=N$ to obtain the existence of $q>0$ such that $$|y|\leq q\implies \sup_{k\in\cK^{t'}}\dbP^{t',k}(\{\Delta^{t',y,\a}_N>\d\} \cup (\dbC^{t'}_{K_\e,1/3})^c)\leq \e.$$
Notice that for all $\o\notin \{\Delta^{t,x,\a}_N>\d\} \cup (\dbC^t_{K_\e,1/3})^c$ one has that 
$$
|\ch_i^{t',y,\a}(\o)-\ch_i^{t',0,\a}|\leq \d,\mbox{ and }|\o_{\ch_i^{t',y,\a}(\o)}-\o_{\ch_i^{t',0,\a}}|\leq K_\e\d^{1/3}.
$$
Thus by the choice of $\d$, if $\ch_N^{t',0,\a}(\o)\wedge \ch_N^{t',0,\a}(\o) =T$ then 
$$
||\hat \o^{t',y,\a}-\hat \o^{t',0,\a}||^{t'}_T\leq K_\e\d^{1/3}\leq \e.
$$
Injecting this to \reff{estimate-finall}, we obtain that for all $\e>0$ there exists $q>0$ such that $$|y|\leq q\implies |\G^{\a}_{n+1} (\pi_n^{(t,x)};t',y) -\G^{\a}_{n+1} (\pi_n^{(t',x)};t',0)|\leq \frac{C\a |t-t'|}{t-t_n}+C \e.$$
Injecting this to (\ref{estimate_gamma_prime}) we obtain that for all $\e>0$ there exists $q>0$ such that
\beaa
|\G^{\a}_{n+1} (\pi_n^{(t,x)};t,0) -\G^{\a}_{n+1} (\pi_n^{(t',x)};t',0)|
\leq C\e +C\sup_{k\in\cK^t}\dbP^{t,k}\left(|B^t_{t'}|\geq q \right)+\rho(\sqrt{|t-t'|}) +C(\e +\a\frac{|t'-t|}{t-t_n})
\eeaa
We now choose $|t-t'|$ small enough to obtain  
$$
|\G^{\a}_{n+1} (\pi_n^{(t,x)};t,0) -\G^{\a}_{n+1} (\pi_n^{(t',x)};t',0)|\leq C\e
$$
which is the continuity we wanted.

\end{document}